\documentclass[a4paper,12pt]{amsart}
\usepackage{amssymb}
\usepackage{ifthen}
\usepackage{graphicx}
\usepackage{amsfonts}
\usepackage{amscd}
\usepackage{amsxtra}
\usepackage{color}

\newtheorem{thm}{Theorem}
\newtheorem{cor}{Corollary}
\newtheorem{lemma}{Lemma}

\newtheorem{conj}{Conjecture}
\theoremstyle{definition}
\newtheorem{example}{Example}%[section]
\newtheorem{prob}{Problem}

\newcounter {own}
\def\theown {\thesection       .\arabic{own}}

\newenvironment{rem}{%
\bigskip
\noindent \textsl{{\sl Remark. }}}{\bigskip}
\newenvironment{rems}{%
\bigskip
\noindent \textsl{{\sl Remarks. }}}{\bigskip}

\newenvironment{pf}[1][]{%
 \vskip 3mm
 \noindent
 \ifthenelse{\equal{#1}{}}%
  {{\slshape Proof. }}%
  {{\slshape #1.} }%
 }%
{\qed\bigskip}

\newcounter{alphabet}
\newcounter{tmp}

\makeatletter
\newcommand{\Ref}[1]{\@ifundefined{r@#1}{}{\setcounter{tmp}{\ref{#1}}\Alph{tmp}}}
\makeatother
\newcommand{\IR}{{\mathbb R}}

\newcommand{\IC}{{\mathbb C}}
\newcommand{\ID}{{\mathbb D}}

\def\be{\begin{equation}}
\def\ee{\end{equation}}

\newcommand{\bee}{\begin{enumerate}}
\newcommand{\eee}{\end{enumerate}}

\newcommand{\blem}{\begin{lemma}}
\newcommand{\elem}{\end{lemma}}
\newcommand{\bthm}{\begin{thm}}
\newcommand{\ethm}{\end{thm}}
\newcommand{\bcor}{\begin{cor}}
\newcommand{\ecor}{\end{cor}}
\newcommand{\beg}{\begin{example}}
\newcommand{\eeg}{\end{example}}
\newcommand{\begs}{\begin{examples}}
\newcommand{\eegs}{\end{examples}}
\newcommand{\bdefe}{\begin{defin}}
\newcommand{\edefe}{\end{defin}}
\newcommand{\bprob}{\begin{prob}}
\newcommand{\eprob}{\end{prob}}
\newcommand{\bei}{\begin{itemize}}
\newcommand{\eei}{\end{itemize}}

\newcommand{\bcon}{\begin{conj}}
\newcommand{\econ}{\end{conj}}
\newcommand{\bcons}{\begin{conjs}}
\newcommand{\econs}{\end{conjs}}
\newcommand{\bprop}{\begin{propo}}
\newcommand{\eprop}{\end{propo}}
\newcommand{\br}{\begin{rem}}
\newcommand{\er}{\end{rem}}
\newcommand{\brs}{\begin{rems}}
\newcommand{\ers}{\end{rems}}
\newcommand{\bo}{\begin{obser}}
\newcommand{\eo}{\end{obser}}
\newcommand{\bos}{\begin{obsers}}
\newcommand{\eos}{\end{obsers}}
\newcommand{\bpf}{\begin{pf}}
\newcommand{\epf}{\end{pf}}
\newcommand{\ba}{\begin{array}}
\newcommand{\ea}{\end{array}}
\newcommand{\beq}{\begin{eqnarray}}
\newcommand{\beqq}{\begin{eqnarray*}}
\newcommand{\eeq}{\end{eqnarray}}
\newcommand{\eeqq}{\end{eqnarray*}}

\def\cc{\setcounter{equation}{0}   % THIS CLEARS THE COUNTER
\setcounter{figure}{0}\setcounter{table}{0}}

\newcounter{minutes}
\divide\time by 60
\newcounter{hours}
\multiply\time by 60 \addtocounter{minutes}{-\time}
\begin{document}
\bibliographystyle{amsplain}

\title[Geometric studies on the class ${\mathcal U}(\lambda)$]{Geometric studies on the class ${\mathcal U}(\lambda)$}
%\title[On the Second coefficient and transformations of analytic functions in the class ${\mathcal U}(\lambda)$]
%{On the Second coefficient and transformations of analytic functions in the class ${\mathcal U}(\lambda)$}

%=========================================================================
\thanks{%$^\dagger$
File:~\jobname .tex,
          printed: \number\day-\number\month-\number\year,
          \thehours.\ifnum\theminutes<10{0}\fi\theminutes}
%=========================================================================

\author[M. Obradovi\'{c}]{Milutin Obradovi\'{c}}
\address{M. Obradovi\'{c},
Department of Mathematics,
Faculty of Civil Engineering, University of Belgrade,
Bulevar Kralja Aleksandra 73, 11000
Belgrade, Serbia }
\email{obrad@grf.bg.ac.rs}

\author[S. Ponnusamy]{Saminathan Ponnusamy $^\dagger $
%${}^{~\mathbf{*}}$
}
\address{S. Ponnusamy,
Indian Statistical Institute (ISI), Chennai Centre, SETS (Society
for Electronic Transactions and Security), MGR Knowledge City, CIT
Campus, Taramani, Chennai 600 113, India }
\email{samy@isichennai.res.in, samy@iitm.ac.in}

 \author[K.-J. Wirths]{Karl-Joachim Wirths}
\address{K.-J. Wirths, Institut f\"ur Analysis und Algebra, TU Braunschweig,
38106 Braunschweig, Germany}
\email{kjwirths@tu-bs.de}

%\thanks{The work of the first author was supported by MNZZS Grant, No. ON174017, Serbia.}
\keywords{Analytic, univalent and starlike functions, Coefficient estimates, subordination, Schwarz' lemma, radius problem,
square-root and $n$-th root transformation. \\
$%{}^{\mathbf{*}}
^\dagger$ {\tt  Corresponding author.
%This author is on leave from Indian Institute of Technology Madras,  India
}
}

\begin{abstract}
The article deals with the family ${\mathcal U}(\lambda)$
of all functions $f$ normalized and analytic in the unit disk such that
$\big |\big (z/f(z)\big )^{2}f'(z)-1\big |<\lambda $  for some $0<\lambda \leq 1$. The family ${\mathcal U}(\lambda)$ has been
studied extensively in the recent past and functions in this family are known to be univalent in $\ID$. However, the problem of determining
sharp bounds for the second coefficients of functions in this family was solved recently in \cite{VY2013} by Vasudevarao and Yanagihara but the proof was
complicated. In this article, we first present a simpler proof. We obtain a number of new subordination results for this family and their consequences.
In addition,  we show that the family  ${\mathcal U}(\lambda )$ is preserved under a number of elementary transformations such as rotation, conjugation, dilation and omitted value transformations, but surprisingly this family is not preserved under the $n$-th root transformation for any $n\geq 2$.  This is a basic here
which helps to generate a number of new theorems and in particular provides a way for constructions of functions from
the family ${\mathcal U}(\lambda)$. Finally, we deal with a radius problem.
\end{abstract}

\maketitle \pagestyle{myheadings}
\markboth{M. Obradovi\'{c}, S. Ponnusamy and K.-J. Wirths}{Geometric studies on the class ${\mathcal U}(\lambda )$}

\cc
\section{Introduction and Basic Properties\label{sec1}}
Let ${\mathcal A}$ be the family of all functions
$f$ analytic in the open unit disk $\ID=\{z\in \IC:\, |z| < 1\}$ with the Taylor series expansion $f(z)=z+\sum_{k=2}^{\infty}a_kz^k$. Let ${\mathcal S}$ denote the subset of ${\mathcal A}$ consisting of functions that are univalent in $\ID$. See \cite{Duren:univ,Go} for the general theory
of univalent functions.  For $0<\lambda \leq 1$, consider the class
$${\mathcal U}(\lambda)=\{f\in {\mathcal A}:\,  \mbox{ $|U_f(z)| <\lambda$ in $\ID$}\},
$$
where $U_f(z)=\big (z/f(z)\big )^{2}f'(z)-1$. Set  ${\mathcal U}:={\mathcal U}(1)$,
${\mathcal U}_2(\lambda):={\mathcal U}(\lambda )\cap \{f\in  {\mathcal A}:\, f''(0)=0\}$ and ${\mathcal U}_2:={\mathcal U}_2(1)$. Because $f'(z)(z/f(z))^2$  $(f\in {\mathcal U})$ is bounded,
it follows that $(z/f(z))^{2}f'(z)\neq 0$ in $\ID$ and thus,
each $f\in {\mathcal U}$ is  non-vanishing in $\ID\backslash\{0\}$.
It is well recognized that the set $\Sigma$ of meromorphic and univalent functions $F$
on  $\{\zeta:\, 1<|\zeta|<\infty \}$ of the form $F(\zeta)=\zeta+\sum_{n=1}^\infty b_n\zeta^{-n}$
plays an indispensable role in the study of ${\mathcal S}$. For $f(z)=1/F(1/z),~\zeta=1/z,$  we have the formula
$$F'(\zeta)=\left(\frac{z}{f(z)}\right)^2f'(z)
$$
and thus, functions in ${\mathcal U}$ are associated with functions $F$ in $\Sigma$ such that $|F'(\zeta)-1|<1$ for $|\zeta|>1$.
In \cite{Aks58}, it was shown that ${\mathcal U}\subsetneqq {\mathcal S}$ and
hence functions in ${\mathcal U}(\lambda)$, that are generalizations of $\mathcal U$, are univalent in $\ID$ for
$0<\lambda \leq 1$. Moreover, if $f\in {\mathcal S}$ and $1/f$ is
a concave schlicht function with the pole at the origin, then $f\in {\mathcal U}$ and this fact is indicated by Aksent\'ev and Avhadiev in
\cite{AksAvh70}. It follows \cite{FR-2006,OP-01,PV2005} that neither ${\mathcal U}$ is
included in ${\mathcal S}^{\star}$ nor includes ${\mathcal S}^{\star}$. Here ${\mathcal S}^{\star}$ denotes the class of starlike functions, namely, functions $f\in{\mathcal S}$ such that $f(\ID)$ is starlike with respect to the origin. In 1995, among many results for the class $\mathcal U$,
Obradovi\'c \cite{Ob-95} proved that if $f\in {\mathcal U}$ then one has the subordination result
$$\frac{z}{f(z)} \prec (1+z)^2, ~  z\in\ID.
$$
For the definition of subordination, denoted by the symbol $\prec$,  we refer to \cite{Duren:univ,Go}.

The class ${\mathcal U}(\lambda)$ has found
many interesting properties \cite{OP-01,obpo-2005, obpo-2007a,obpo-2009a,ObPo-Jammu,VY2013}.
It is a simple exercise to see that
each $f\in{\mathcal U}(\lambda)$ has the characterization \cite{obpo-2007a}
\be\label{OPW7-eq1}
\frac{z}{f(z)}=1-a_{2}z+\lambda z\int_{0}^{z}\omega(t)\,dt,
\ee
for some  $\omega\in{\mathcal B}$, where $a_2 =f''(0)/2$, and ${\mathcal B}$ denotes the
class of functions $\omega$ analytic in $\ID$ such that $|\omega(z)|\leq 1$ for $z\in\ID$.
%$$U_f(z)=\big (z/f(z)\big )^{2}f'(z)-1.
%$$
%
Here is a typical set of functions in $\mathcal{U}\cap \mathcal{S}^*$ given by
$$L=\left \{z, ~~\frac{z}{(1\pm z)^2},~~\frac{z}{1\pm z},~~\frac{z}{1\pm z^2},
~~ \frac{z}{1\pm z+z^2}\right \},
$$
where $L$ is exactly the set of functions in $\mathcal{S}$ having integral coefficients in the power series expansion, \cite{Fr1}.
Since ${\mathcal U}\subsetneq {\mathcal S}$ and the Koebe function $z/(1-z)^2$ belongs to ${\mathcal U}$,
$|a_2|\leq 2$ is obvious for $f\in {\mathcal U}$. The sharp estimation for the second coefficient of
functions in ${\mathcal U}(\lambda)$ was known only recently in \cite{VY2013}. One of our main aims in this article
is to give a simpler and different proof of this result. More precisely, in Theorem \ref{OPW7-th2}, we present a new proof that
if ${\mathcal U}(\lambda)$, then $|a_2|\leq 1+\lambda$ holds, and,
in Theorem \ref{OPW7-th3}, we show that if $|a_2|=1+\lambda$, then $f$ must be of the form
\be\label{ex1}
f(z)= \frac{z}{1-a_2z+\lambda  e^{i\theta} z^2}
\ee
for some $\theta \in [0,2\pi]$.

It is well-known that the class $\mathcal S$ is preserved under a number of elementary
transformations, eg. conjugation, rotation, dilation, disk automorphisms (i.e. the Koebe transformations),
range, omitted-value and square-root transformations to say few.
% We show that ${\mathcal U}(\lambda )$ as a subset of $\mathcal S$
%preserves some of these properties and as a consequence we derive few applications.

\blem \label{lem1}
The class ${\mathcal U}(\lambda )$ is preserved under rotation, conjugation, dilation and omitted-value transformations.
\elem
\bpf Let $f\in {\mathcal U}(\lambda )$ and define $g(z)=e^{-i\theta}f(ze^{i\theta})$, $h(z)=\overline{f(\overline{z}})$
and $\psi(z)=r^{-1}f(rz)$. Then we see that
$g'(z)=f'(ze^{i\theta})$, $h'(z)=\overline{f'(\overline{z}})$, $\psi '(z)=f'(rz)$,
\beqq
\left (\frac{z}{g(z)} \right )^{2}g'(z) -1 &= &\left (\frac{ze^{i\theta}}{f(ze^{i\theta})} \right )^{2}f'(ze^{i\theta}) -1,\\
\left (\frac{z}{h(z)} \right )^{2}h'(z) -1 &=&
%\left (\frac{z}{\overline{f(\overline{z}})} \right )^{2}\overline{f'(\overline{z})} -1=
\overline{\left (\frac{\overline{z}}{f(\overline{z})} \right )^{2}f'(\overline{z})}-1, ~\mbox{ and}\\
\left (\frac{z}{\psi(z)} \right )^{2}\psi'(z) -1 &= &\left (\frac{rz}{f(rz)} \right )^{2}f'(rz) -1.
\eeqq
It follows that $g$, $h$ and $\psi$ belong  to ${\mathcal U}(\lambda )$.

Finally, if $f\in {\mathcal U}(\lambda )$ and $f(z)\neq c$ for some $c\neq 0$, then the function $F$ defined by
$$F(z)=\frac{cf(z)}{c-f(z)}
$$
obviously belongs to $\mathcal S$. Thus, $z/F(z)$ is non-vanishing in $\ID$, and it is a simple exercise to see that
\be\label{OS-eq2}
U_f(z)=\left (\frac{z}{f(z)} \right )^{2}f'(z)-1=\frac{z}{f(z)} -z\left (\frac{z}{f(z)} \right )'-1, \quad z\in\ID.
\ee
Using \eqref{OS-eq2}, it is easy to see that $U_{F}(z)=U_{f}(z)$ for $z\in \ID$. Consequently, $F\in {\mathcal U}(\lambda )$.
The proof is complete.
\epf

\bcor\label{13-th6}
Let $f\in {\mathcal U}(\lambda)$ for some $0<\lambda \leq1$  and $a_2 =f''(0)/2$. If $a_{2}+\mu\neq 0$ for some complex number
$\mu$ with $|\mu|\leq 1-\lambda$, then
$$ -\frac{1}{a_{2}+\mu}\notin f(\ID).
$$
\ecor
\bpf
Let $f\in {\mathcal U}(\lambda)$. Suppose that there exists a point $z_{0}\in\ID$ such that $f(z_{0})=- \frac{1}{a_{2}+\mu}$. Then
$$\frac{z_{0}}{f(z_{0})}=- (a_{2}+\mu )z_{0}
$$
and thus, according to the representation \eqref{OPW7-eq1}, the last relation implies that
$$1+\mu z_{0}+\lambda z_{0}\int_{0}^{z_{0}}\omega (t)\,dt=0
$$
for some $\omega\in{\mathcal B}$. %analytic function $\omega$ in $\ID$ with $|\omega(z)|\leq 1$ in $\ID$.
 But, this is not possible because
\beqq
\left|1+\mu z_{0}+\lambda z_{0}\int_{0}^{z_{0}}\omega (t)\,dt \right
|&\geq& 1-|\mu|\,|z_{0}|-\lambda |z_{0}|^{2}\\
&\geq& 1-(1-\lambda)|z_0|-\lambda |z_{0}|^{2}\\
&=& (1-|z_{0}|)(1+\lambda |z_{0}|)>0.
\eeqq
We complete the proof.
\epf

According to Corollary \ref{13-th6}, the function $F$ defined by
$$ F(z)= \frac{f(z)}{1+(a_{2}+\mu)f(z)}
$$
belongs to the class ${\mathcal U}(\lambda)$ whenever $f\in {\mathcal U}(\lambda)$ and $a_{2}+\mu\neq 0$ with $|\mu|\leq 1-\lambda$. In particular,
$$ F(z)= \frac{f(z)}{1+(a_{2}+1-\lambda )f(z)}
$$
belongs to the class ${\mathcal U}(\lambda)$ if $f\in {\mathcal U}(\lambda)$ and $a_{2}\neq \lambda -1$.

On the other hand,  the class ${\mathcal U}$ (and hence, ${\mathcal U}(\lambda )$)
is not preserved under the square-root transformation. For example, we consider the function
$$f_1(z)=\frac{z}{1+(1/2)z+(1/3)z^3}.
$$
Then we see that $z/f_1(z)$ is non-vanishing in $\ID$, and it is a simple exercise to see that
$U_{f_1}(z)=-(2/3)z^3$ showing that $f_1\in{\mathcal U}$. In particular, $f_1$ is univalent in $\ID$. On the other hand if we consider $g_1$
defined by
$$ g_1(z)=\sqrt{f_1(z^2)}=z\sqrt{\frac{f_1(z^2)}{z^2}}
$$
then, because $\mathcal S$ is preserved under the square-root transformation, it follows that $g_1$ is univalent in $\ID$ whereas
$$\left (\frac{z}{g_1(z)} \right )^{2}g_1'(z)  -1= \left (\frac{z}{f_1(z)} \right )^{3/2}f_1'(z)  -1
=\frac{1-(2/3)z^6}{\sqrt{1+(1/2)z^2+(1/3)z^6}} -1
$$
which approaches the value $\frac{5\sqrt{6}-3}{3} >1$ as $z\rightarrow i$. This means that $U_{g_1}(\ID)$ cannot be a subset of the unit disk $\ID$ and hence,
the square-root transformation $g_1$ of $f_1$ does not belong to ${\mathcal U}$.

More generally if we consider
$$f(z)=\frac{z}{1+(1/n)z+(-1)^n(1/(n+1))z^{n+1}}
$$
then a computation shows that $f\in{\mathcal U}$ whereas the $n$-th root transformation $g$ of $f$, given by
$$g(z)=\sqrt[n]{f(z^n)}=z\sqrt[n]{\frac{f(z^n)}{z^n}},
$$
does not belong to the class ${\mathcal U}$ for each $n\geq 2$. Thus, for any $n\geq 2$, ${\mathcal U}$ is not preserved under
the $n$-th root transformation unlike the class ${\mathcal S}$.

The remaining part of the article is organized as follows. In Section \ref{sec2}, we present a sharp coefficient bound for the second
Taylor coefficient of  $f\in {\mathcal U}(\lambda)$ and prove, in particular,  several subordination results for $z/f(z)$ implying growth theorems for
the family ${\mathcal U}(\lambda)$.
In Section \ref{sec3}, we derive subordination results for functions in the family ${\mathcal U}(\lambda)$ and in Section \ref{sec4}, we present
a number of consequences of Lemma \ref{lem1}. Section \ref{sec5} is dedicated to examples
of construction principles for functions in ${\mathcal U}(\lambda)$. The aim of Section \ref{sec6} is
the calculation of a radius $r_0$ such that $f(r_0 z)/r_0$ belongs to ${\mathcal U}$ if
$f$ is univalent in the unit disk.

\section{Second Coefficient for functions in  ${\mathcal U}(\lambda)$\label{sec2}}

First we present a direct approach and later we shall obtain the following result as a simple consequence of a subordination
result (see Theorems \ref{OPW7-th4} and \ref{OPW7-th2a}).

\bthm\label{OPW7-th2}
Let $f\in {\mathcal U}(\lambda)$ for some $0<\lambda \leq 1$.  Then $|a_2|\leq 1+\lambda$.
\ethm
\bpf
Recall the fact that $f(z)=z+\sum_{n=2}^{\infty}a_nz^n \in {\mathcal U}(\lambda)$ if and only if
\be\label{OPW7-eq3}
 \frac{z}{f(z)}=1-a_2z+\lambda z\int_0^z \omega(t)\,dt \,\neq 0, \quad z\in \ID,
\ee
where $\omega\in{\mathcal B}$.

It suffices to prove that for $|a_2|>1+\lambda$ and for any $\omega \in{\mathcal B}$, there exists a $z_0\in \ID$ such that
\[ 1-a_2z_0+\lambda z_0\int_0^{z_0} \omega(t)\,dt= 0.
\]

We may now assume that
\be\label{OPW7-eq4}
|a_2|=\frac{1+\lambda}{r},\quad r\in (0,1),
\ee
and prove that the map $F$ defined by
\[ a_2F(z)=1+\lambda z\int_0^z \omega(t)\,dt
\]
is a contracting map of $\ID_r$ into $\ID_r$, where $\ID_r =\{z:\, |z|\leq r\}$.

We see that for $z\in \ID_r$,
\[ |F(z)| =  \frac{r}{1+\lambda} \left |1 + \lambda z\int_0^z \omega(t)\,dt\right |\leq \frac{r(1+\lambda |z|^2)}{1+\lambda} < r.
\]
Now let $z_1,z_2 \in  D_r$. This gives that
\beqq
|F(z_1)-F(z_2)| &= &  \frac{\lambda r}{1+\lambda}\left|z_1\int_0^{z_1} \omega(t)\,dt  +(-z_1+z_1-z_2)\int_0^{z_2} \omega(t)\,dt\right|\\
&\leq & \frac{\lambda r}{1+\lambda}\left(|z_1|\left|\int_{z_2}^{z_1}\omega(t)\,dt\right| + |z_1-z_2|\left|\int_0^{z_2}\omega(t)\,dt\right|\right)\\
&\leq & \frac{\lambda r}{1+\lambda}\,(|z_1|+|z_2|)|z_1-z_2|\\
& \leq& r^2 |z_1-z_2|.
\eeqq
Thus,  $F$ is a contracting map of $\ID_r$ into $\ID_r$. This implies, according to Banach's fixed point theorem, that there exists a
$z_0\in \ID_r$ such that $F(z_0)=z_0$ which contradicts \eqref{OPW7-eq3} at $z_0\in \ID$ (and thus, \eqref{OPW7-eq4} is not true for any $r\in (0,1)$). Hence, we must have
$|a_2|\leq 1+\lambda$ for $f\in {\mathcal U}(\lambda)$.
\epf

Determining the sharp bound for the Taylor coefficients $|a_n|$ ($n\ge 3$), for $f\in {\mathcal U}(\lambda)$, remains an open problem.

Next we deal with the equality case.

\bthm\label{OPW7-th3}
If  $f\in {\mathcal U}(\lambda)$ and   $|a_2|=1+\lambda$, then $f$ must be of the form \eqref{ex1} and especially,
$$f(z)=\frac{z}{1-(1+\lambda)e^{i\phi}z+\lambda e^{2i\phi}z^2}.
$$
%$$f(z)= \frac{z}{1-a_2z+\lambda  e^{i\theta} z^2}
%$$
 %and
%\[
% \frac{z}{f(z)}=1-a_2z+\lambda z\int_0^z \omega(t)\,dt, \quad z\in \mathbb{D},
%\]
%for some $\theta \in [0,2\pi]$.
\ethm\bpf
Let $f\in {\mathcal U}(\lambda)$. Then $f$ must be of the form \eqref{OPW7-eq3}
% \frac{z}{f(z)}=1-a_2z+\lambda z\int_0^z \omega(t)\,dt \neq 0, \quad z\in \ID,
for some $\omega\in{\mathcal B}$.
If $|a_2|=1+\lambda$, then we must show that $\omega$ in \eqref{OPW7-eq3} takes the form $\omega(z)= e^{i\theta}$ for some $\theta \in [0,2\pi]$ and all  $z\in \mathbb{D}.$

Assume on the contrary that $\omega (0)=a \in \mathbb{D}$ and $f$ as in \eqref{OPW7-eq3}. Then, according to Schwarz-Pick's Lemma applied to $\omega\in{\mathcal B}$, we get
\[ \left|\frac{a-\omega(z)}{1-\overline{a}\omega(z)}\right|\leq |z|,\quad z\in \mathbb{D},
\]
%which is equivalent to the inequality
%\[
% \left|\omega(z)-a\frac{1-|z|^2}{1-|a|^2|z|^2}\right|\leq |z|\frac{1-|a|^2}{1-|a|^2|z|^2},\quad z\in \mathbb{D}.
%\]
%From this inequality, we get immediately
from which we can immediately obtain that
\[
 |\omega(z)|\leq \frac{|a|+|z|}{1+|az|},\quad z\in \mathbb{D},
\]
and thus, we see that
\beqq
\left|\int_0^z \omega(t)\,dt\right| &\leq & \int_0^{|z|}\frac{|a|+s}{1+|a|s}\,ds
\, = \frac{|z|}{|a|}-\frac{1-|a|^2}{|a|^2}\log(1+|az|)\\
&\leq& \frac{1}{|a|}-\frac{1-|a|^2}{|a|^2}\log(1+|a|)=:v(|a|)<1.
\eeqq
Now, we let as in Theorem \ref{OPW7-th2},
\[
 F(z)=\frac{1+\lambda z\int_0^z \omega(t)\,dt}{a_2}
\]
and define
\[
 \frac{1+\lambda v(|a|)}{1+\lambda}=:r<1.
\]
For $z\in \ID_r$ we have
\[
 |F(z)|\leq \frac{1+\lambda r v(|a|)}{1+\lambda}<r,
\]
and for $z_1,z_2\in \ID_r$ we get as above
\beqq
|F(z_1)-F(z_2)|&=&\frac{\lambda}{1+\lambda}\left|z_1\int_{z_2}^{z_1}\omega(t)\,dt+(z_1-z_2)\int_0^{z_2}\omega(t)\,dt\right|\\
&\leq & \frac{1}{2}(|z_1|+|z_2|)|z_1-z_2|\leq r |z_1-z_2|.
\eeqq
Hence $F$ has a fixed point in $\ID_r$ which contradicts $f \in {\mathcal U}(\lambda)$.

At last, we consider for fixed $\varphi, \psi \in [0, 2\pi]$  the cases
\[
 \frac{z}{f(z)}=1-(1+\lambda)e^{i\varphi}z+\lambda e^{i\psi}z^2=:p(\varphi, \psi, z)
\]
and prove that $p(\varphi, \psi, z) $ is nonvanishing in the unit disk if and only if $\psi = 2\varphi$.

Without restriction of generality we may assume $\varphi = 0$ and prove that among the functions $p(0, \psi, z) $ the only one
non-vanishing in $\mathbb{D}$ is the function $p(0,0,z)$.

To that end we consider the functions
\[
 q_{\psi}(z):=(1+\lambda)z-\lambda e^{i\psi}z^2.
\]
Since for $z=re^{i\tau}$, $r\in [0,1),\, \tau \in [0,2\pi],$ the inequality
\[
 {\rm Re}\, q'_{\psi}(z)=1+\lambda-2\lambda r \cos(\psi+\tau) >0
\]
is valid, the function $q_{\psi}$ is univalent in $\mathbb{D}$. In our case $q_{\psi}(\partial \mathbb{D})$ is a Jordan curve and $q_{\psi}( \mathbb{D})$ is the simply connected domain bounded by this curve. If we consider the curve $q_{\psi}(\partial \mathbb{D})$, we see that
\[
 \left|q_{\psi}\big(e^{i\tau}\big)\right|\geq 1+\lambda-\lambda=1, \quad \tau \in [0,2\pi],
\]
and the minimum modulus is attained if and only if $e^{i\tau}=e^{i(\psi + 2\tau)}$, i.e. $\tau=-\psi$. Hence,
$1 \notin q_{\psi}( \mathbb{D})$, if and only if
\[
 {\rm Re}\,q_{\psi}\big(e^{-i\psi}\big)= (1+\lambda)\cos \psi - \lambda \cos \psi =\cos \psi =1.
\]
This is satisfied if and only if $\psi = 0.$ Thus, $f$ must be of the form \eqref{ex1}.
%$$f(z)= \frac{z}{1-a_2z+\lambda  e^{i\theta} z^2}
%$$
%and we complete the proof.
\epf

\section{Subordination \label{sec3}}
\bthm\label{OPW7-th1}
Let $f\in {\mathcal U}(\lambda)$ for some $0<\lambda \leq 1$ and $a_2 =f''(0)/2$. Then
$$ \frac{z}{f(z)}+ a_2z \prec 1+2\lambda z + \lambda z^2.
$$
%In particular, if $f\in {\mathcal U}_2$ then $\frac{z}{f(z)} \prec (1+z)^2$.
\ethm
\bpf
From \eqref{OPW7-eq1}, we observe that each $f\in {\mathcal U}(\lambda)$ has the form
\be\label{OPW7-eq2}
 \frac{z}{f(z)}=1-a_2z+\lambda \psi (z),  \quad \psi (z)=z\int_0^z\,\omega(t)\,dt,
\ee
where $\omega\in{\mathcal B}$. Since $|\omega(z)|\leq 1$ for $z\in\ID$ and $\phi (z)=\psi (z)/z$ has the property that $\phi (0)=0$ and $|\phi (z)|\leq 1$, the classical Schwarz' lemma shows that
%the relation
%$$\psi (z)=z\int_0^z\,\omega(t)\,dt=z^2\int_0^1\,\omega(uz)\,du
%$$
$|\psi (z)|\leq |z|^2$ in $\ID$.
%This observation shows that the function $g$ defined by
%\[ \frac{z}{g(z)}=1+\psi (z)
%\]
%is non-vanishing in $\ID$, and $g\in {\mathcal U}_2$. Because
%$$\frac{z^2}{2}\prec z+ \frac{z^2}{2} ~\mbox{ and }~ \left |\frac{z}{g(z)} -1\right |\leq |z|^2,
%$$
%it follows that  $\psi (z) \prec 2z+z^2$ in $\ID$, i.e. $\frac{z}{g(z)}\prec 1+2z+z^2= (1+z)^2$.
Again, because
$$\frac{z^2}{2}\prec z+ \frac{z^2}{2} ~\mbox{ and }~ |\psi (z)|\leq |z|^2,
$$
it follows that  $\psi (z) \prec 2z+z^2$ in $\ID$. The desired conclusion follows from  \eqref{OPW7-eq2}.
\epf

As remarked earlier, our next result includes a proof of Theorem \ref{OPW7-th2} which will be stated as a corollary below.

\bthm\label{OPW7-th4}
If  $f\in {\mathcal U}(\lambda)$ for $\lambda \in (0,1]$, then
\be\label{OPW7-eq4c}
\frac{f(z)}{z}\prec \frac{1}{1+(1+\lambda)z+\lambda z^2}, ~ z\in \ID ,
\ee
or equivalently
\[ \frac{z}{f(z)}\prec 1+(1+\lambda)z+\lambda z^2, ~ z\in \ID;
\]
and, for $|z|=r$,
$$\left|\frac{z}{f(z)}-1\right|\leq -1+(1+\lambda r)(1+r).
$$
In particular, if $f\in {\mathcal U}$ then $\frac{z}{f(z)} \prec (1+z)^2$ in $\ID$.
\ethm
\bpf It suffices to prove the theorem for $\lambda \in (0,1)$.
Assume that $f\in {\mathcal U}(\lambda)$ and $s(z)=1+(1+\lambda)z+\lambda z^2$. First we observe that $s(z)$ is univalent in $\ID$ for
$\lambda \in (0,1)$. Indeed for $z_1,z_2$ in the closed unit disk $\overline{\mathbb{D}}$, we have
\[\left|\frac{s(z_1)-s(z_2)}{z_1-z_2}\right|=|1+\lambda +\lambda (z_1+z_2)| \geq  1+\lambda - 2\lambda >0\]
(and also ${\rm Re\,}s'(z)\geq 1+\lambda - 2\lambda >0$ in $\overline{\ID}$) showing that $s(z)$ is univalent in $\ID$.

We need to show that $\frac{z}{f(z)}\prec s(z)$. Suppose on the contrary that $\frac{z}{f(z)}$ is not subordinated to $s(z)$.
%To apply Lemma 2.2.d of the book of Miller and Mocanu, we show that $s$ is injective on $\overline{\mathbb{D}}$. This follows from the inequalities
%\[
%\left|\frac{s(z_1)-s(z_2)}{z_1-z_2}\right|=|1+\lambda +\lambda (z_1+z_2)| \geq  1+\lambda - 2\lambda >0,~z_1,z_2\in\overline{\mathbb{D}}.\]
%The application of the above lemma
As an application of \cite[Lemma 1]{MM-MMJ85} (see also \cite{MM-MMJ81}), there exist points $z_0=r_0e^{i\theta_0}\in {\mathbb{D}}$ and $\zeta_0 \in \partial{\mathbb{D}}$ such that
\[\frac{z_0}{f(z_0)}=1+ (1+\lambda)\zeta_0 + \lambda \zeta_0^2.
\]
On the other hand we know from \cite[Theorem 3.2]{VY2013} that $\frac{z_0}{f(z_0)}$ lies in the union of the images of the disks $\{z:\, |z|\leq r_0\}$ under the functions
\begin{equation}\label{f1}
\frac{z}{g(z)}=1 +(1+\lambda e^{i\tau})z +\lambda e^{i\varphi}z^2
\end{equation}
where one has to consider only those $g$ belonging to ${\mathcal U}(\lambda) $. Hence, for our purposes it is sufficient to prove that the functions of
the type (\ref{f1}), where $g$ is restricted as above, are subordinated to the function $s(z)$. We observe that functions of the type $g$ given by (\ref{f1}) belong to ${\mathcal U}(\lambda)$ if and only if
\begin{equation}\label{f2}
0\neq 1+(1+\lambda e^{i\tau})z+\lambda e^{i\varphi}z^2,\quad z\in \mathbb{D}.
\end{equation}
Using the abbreviation
\[1+\lambda e^{i\tau}=|1+\lambda e^{i\tau}|e^{i\gamma}\]
we get
\[(1+\lambda e^{i\tau})z+\lambda e^{i\varphi}z^2=e^{i(2\gamma - \varphi)}(|1+\lambda e^{i\tau}|e^{i(\varphi-\gamma)}z+\lambda e^{2i(\varphi-\gamma)}z^2).\]
Hence, (\ref{f2}) is equivalent to
\[
-e^{-i((2\gamma - \varphi)}\neq |1+\lambda e^{i\tau}| u+\lambda u^2, \quad u \in \mathbb{D}.\]
In the following we let $\beta=\varphi-2\gamma$ and \[l=|1+\lambda e^{i\tau}|\in [1-\lambda,1+\lambda].\]
For $u=e^{i\alpha}$ and $x+iy=le^{i\alpha}+\lambda e^{2i\alpha},$
we have
\begin{equation}\label{f2-a}
x+\lambda=\cos \alpha (l+2\lambda \cos \alpha) ~\mbox{ and }~ y=\sin \alpha (l+2\lambda \cos \alpha).
\end{equation}
\begin{figure}[H]
\begin{center}
\includegraphics[height=6.0cm, width=5.5cm, scale=1]{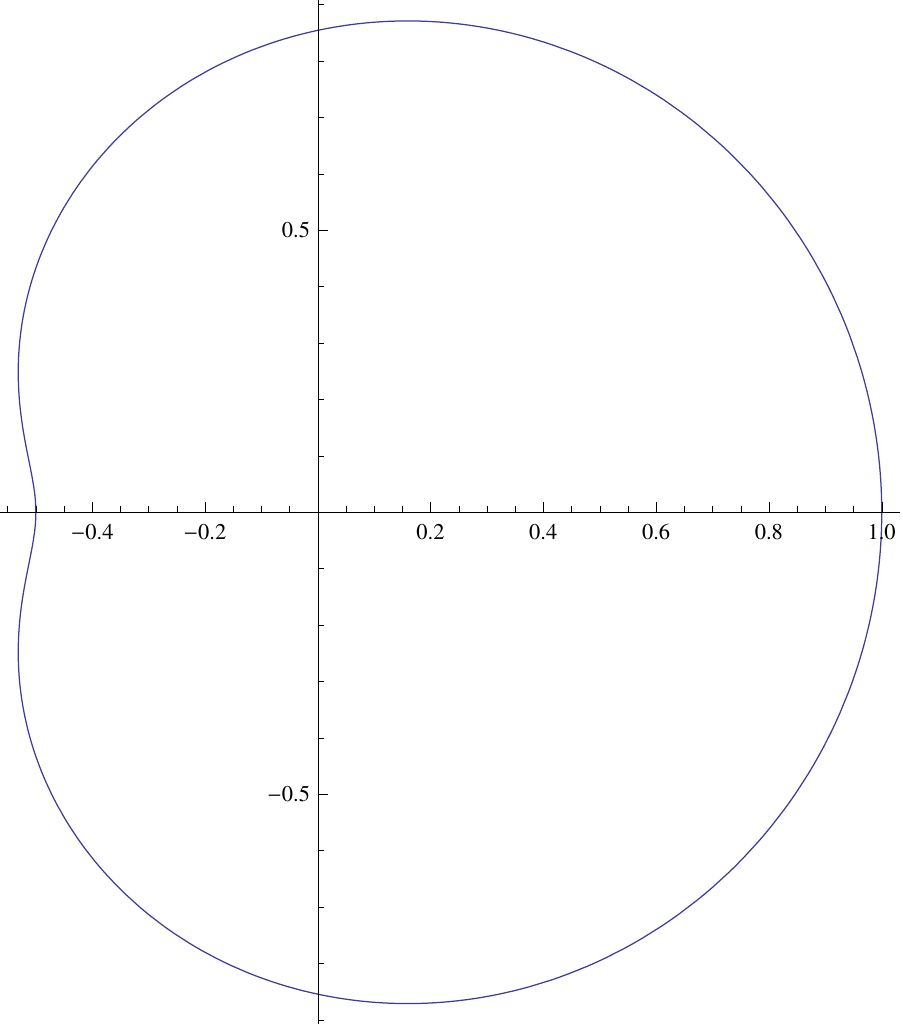}
\hspace{1cm}
\includegraphics[height=6.0cm, width=5.5cm, scale=1]{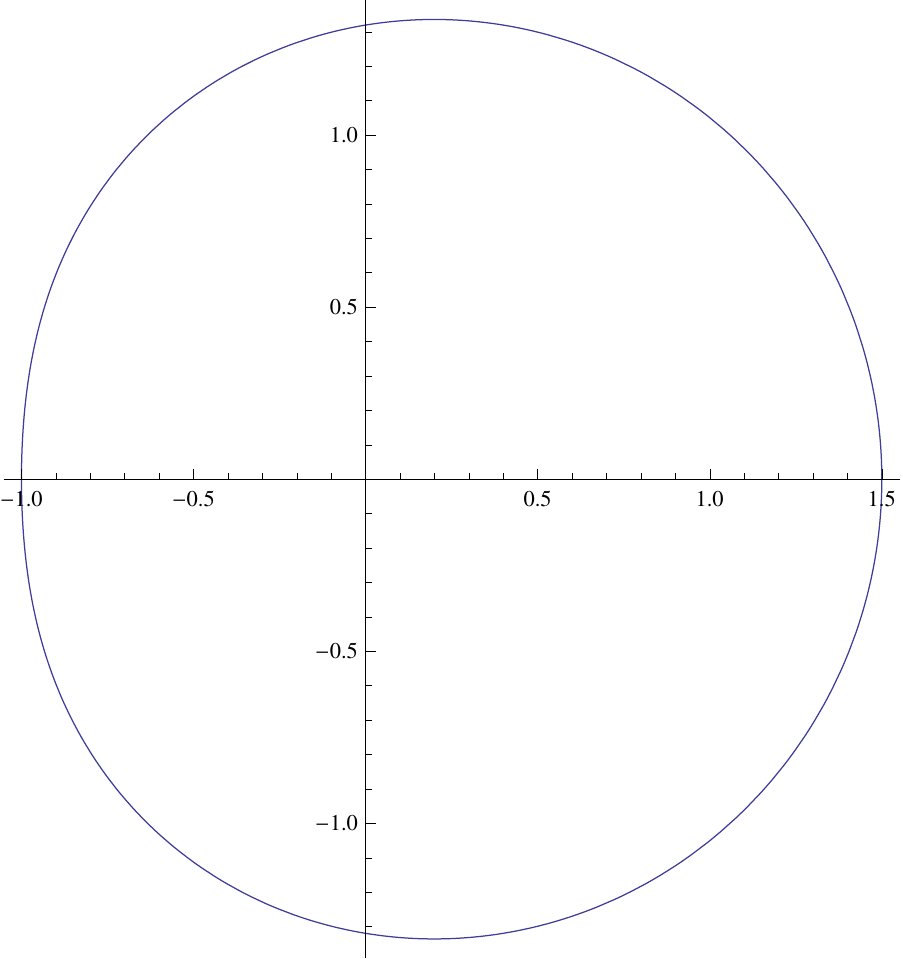}
\end{center}
$\lambda=0.25$, $l=0.75$ \hspace{4cm} $\lambda=0.25$, $l=1.25$

\vspace{0.25cm}

\begin{center}
\includegraphics[height=6.0cm, width=5.5cm, scale=1]{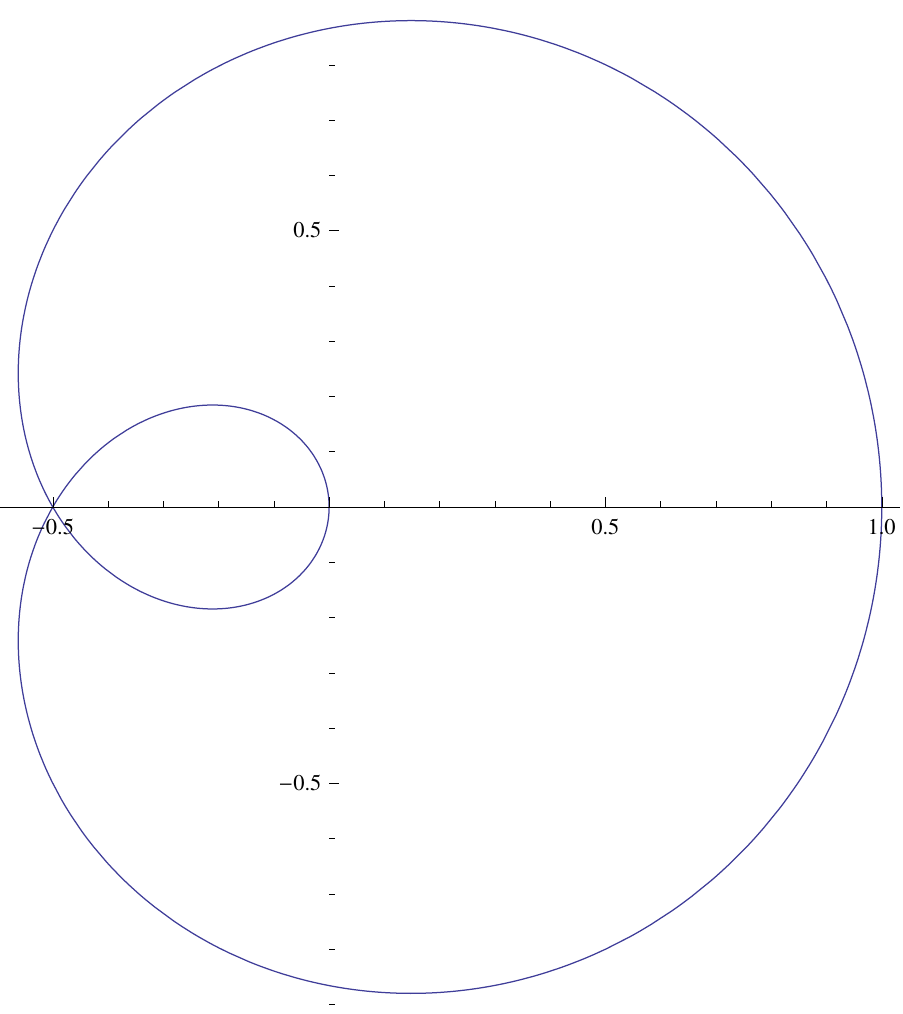}
\hspace{1cm}
\includegraphics[height=6.0cm, width=5.5cm, scale=1]{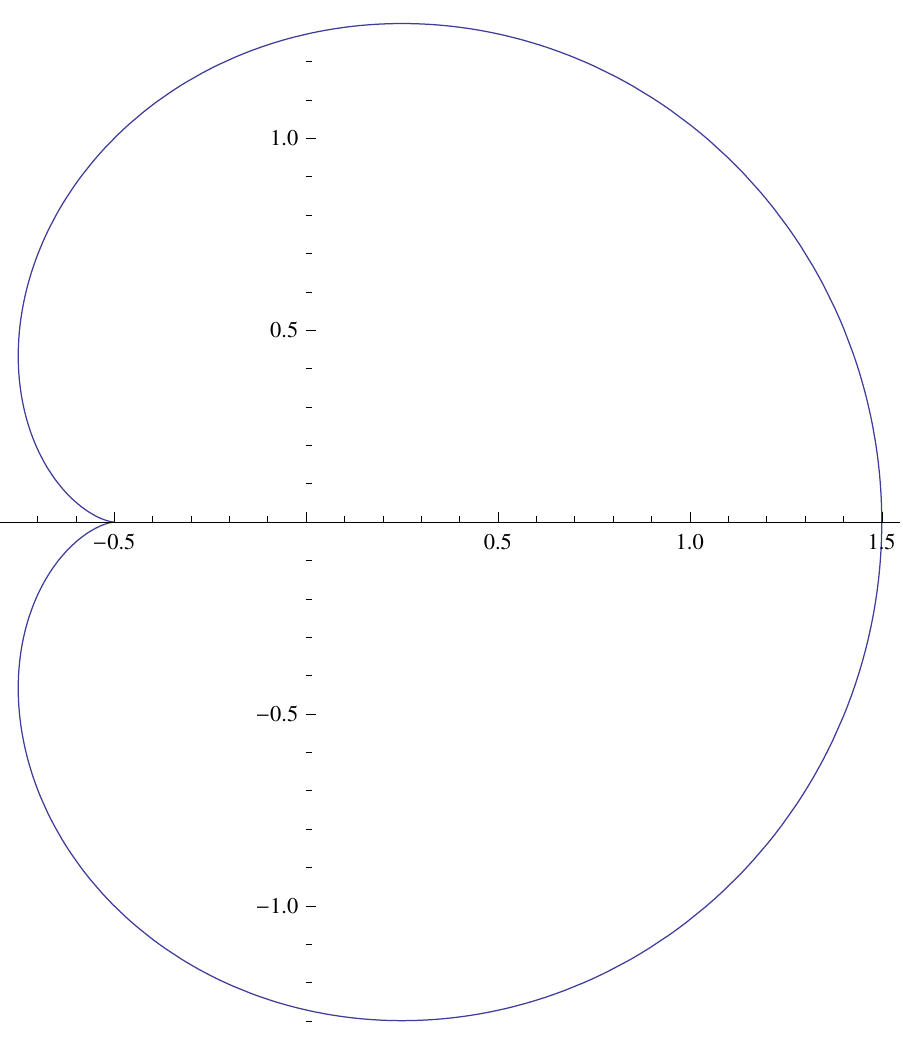}
\end{center}
$\lambda=0.5$, $l=0.5$ \hspace{4cm} $\lambda=0.5$, $l=1$

\vspace{0.25cm}

\begin{center}
\includegraphics[height=6.0cm, width=5.5cm, scale=1]{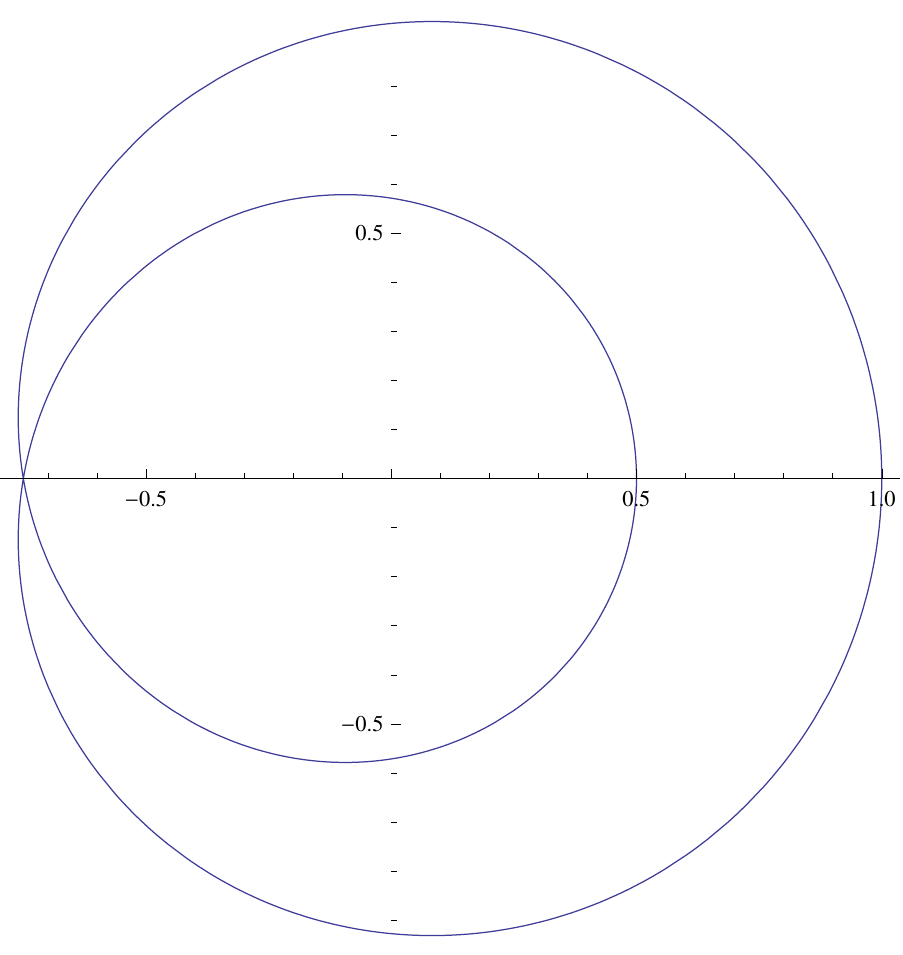}
\hspace{1cm}
\includegraphics[height=6.0cm, width=5.5cm, scale=1]{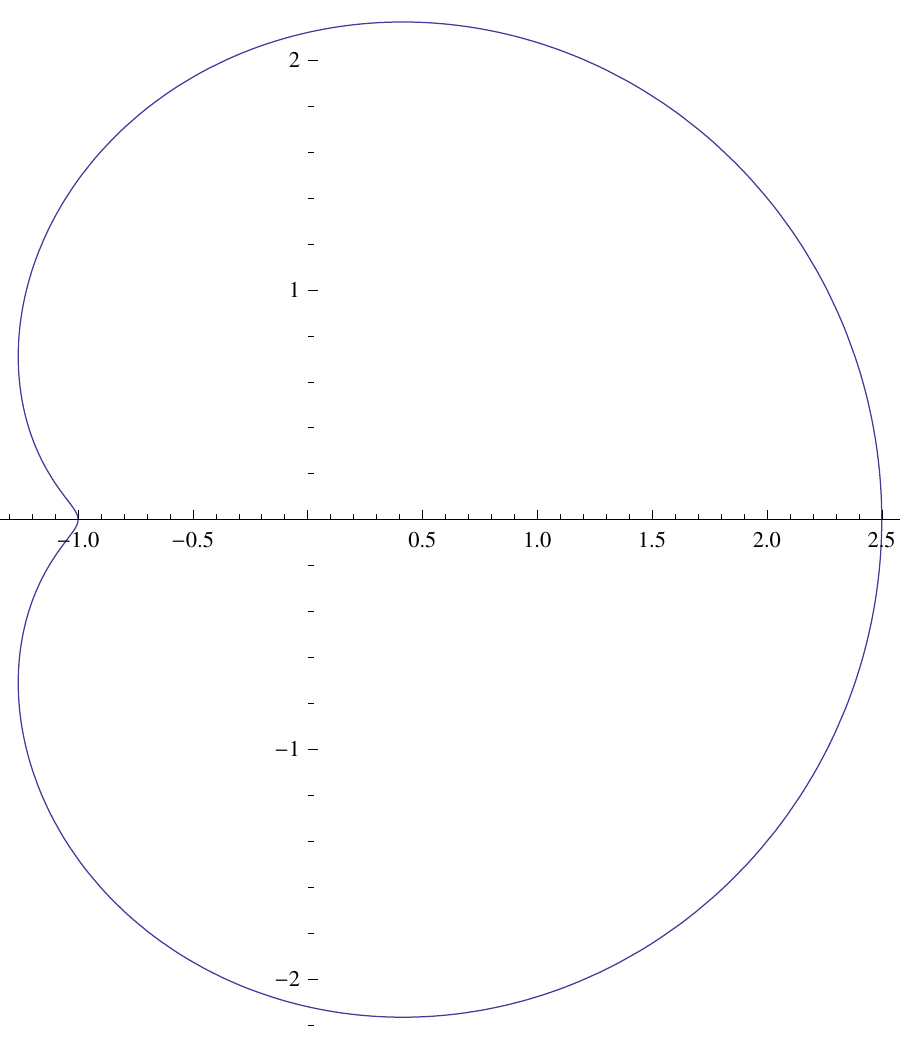}
\end{center}
$\lambda=0.75$, $l=0.25$ \hspace{4cm} $\lambda=0.75$, $l=1.75$
\caption{The graph of some lima\c{c}ons parameterized by  \eqref{f2-a}  for certain values of $\lambda$ and $l$.\label{limacon}}
\end{figure}
This is the parametrization of a lima\c{c}on with center $(-\lambda, 0)$ (see Figure \ref{limacon} for the graph of some
lima\c{c}ons parameterized by \eqref{f2-a} for various values of $\lambda$ and $l$).
The implicit equation of this lima\c{c}on derived from the above equations is the following
\[
(x^2+y^2-\lambda^2)^2=l^2(x^2+y^2+\lambda^2+2\lambda x).\]
The intersection points $(x,y)$ of the lima\c{c}on and the unit circle can be got from this equation and
\[
\frac{(1-\lambda^2)^2-l^2(1+\lambda^2)}{2\lambda l^2}=x=:-\cos \beta_1.\]
Hence, for $|\beta |\leq \beta_1$ the functions $g$ defined by (\ref{f1}) belong to ${\mathcal U}(\lambda)$.

For $l=1+\lambda$, the case $\varphi = 0$ is the only one that produces a member of ${\mathcal U}(\lambda)$ in (\ref{f1}), whereas for $l=1-\lambda$ all functions $g$ defined by (\ref{f1}) belong to this family.

\begin{figure}[H]
\begin{center}
\includegraphics[height=6.0cm, width=5.5cm, scale=1]{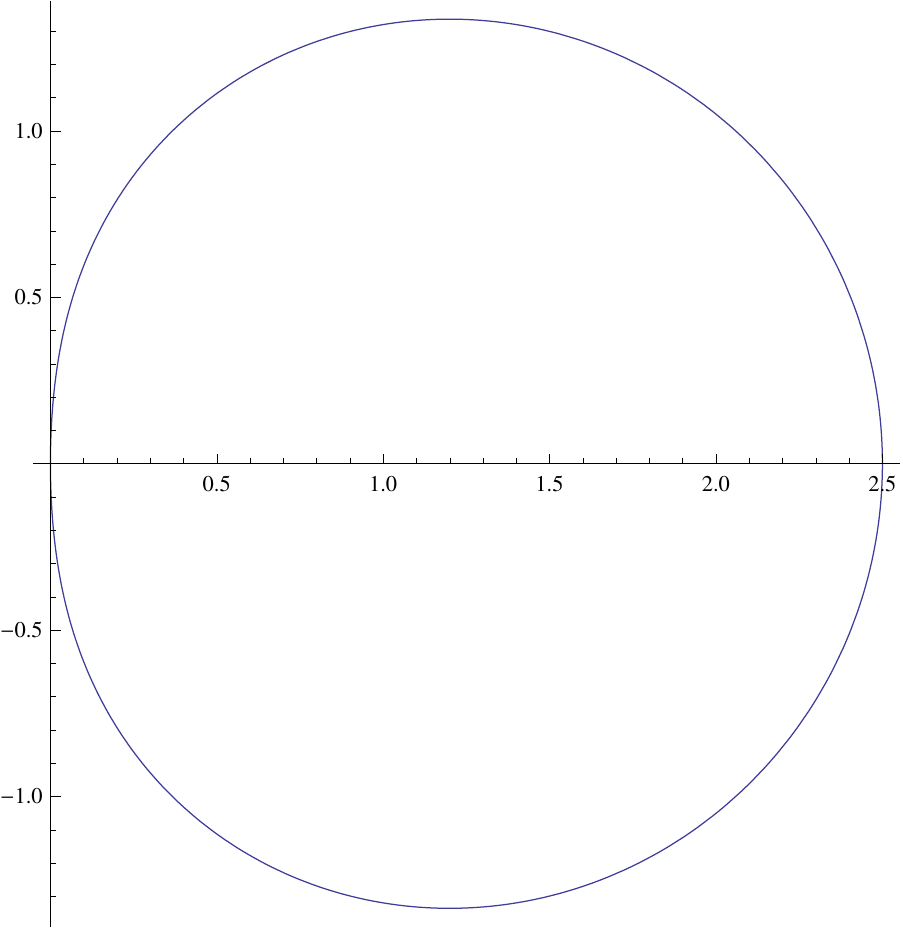}
\hspace{1cm}
\includegraphics[height=6.0cm, width=5.5cm, scale=1]{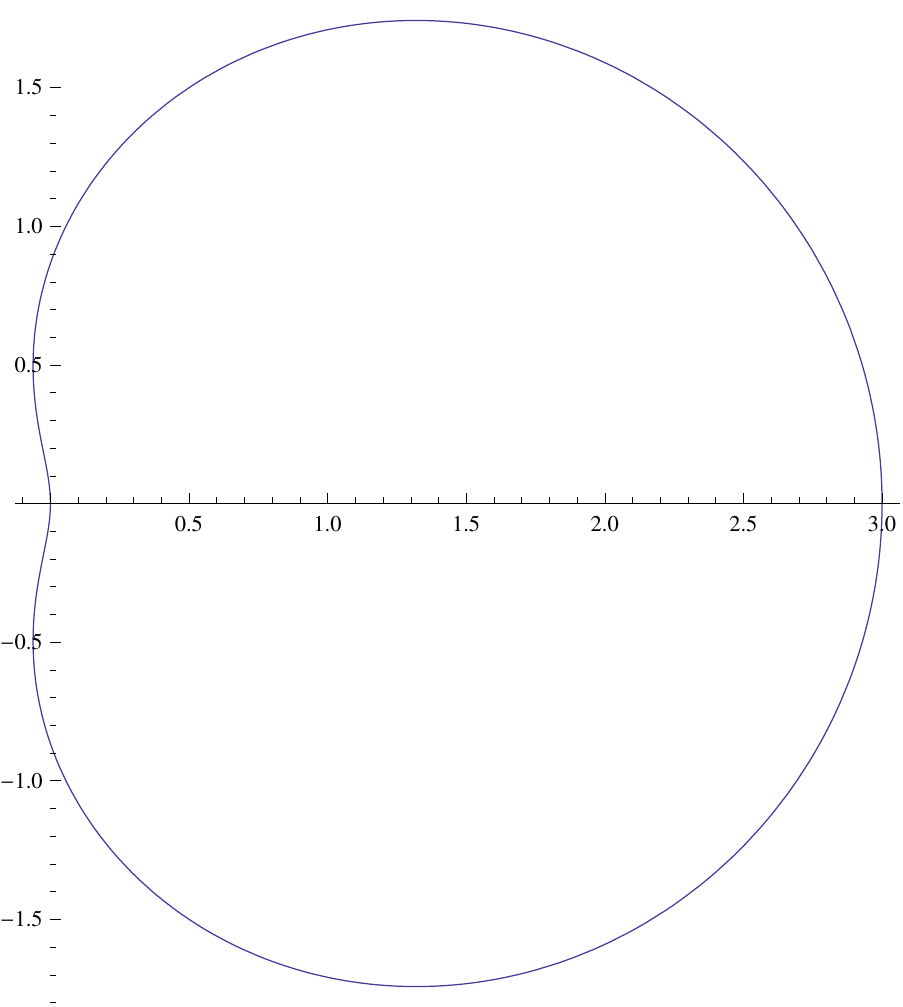}
\end{center}
$\lambda=0.25$ \hspace{5cm} $\lambda=0.5$

\vspace{0.25cm}

\begin{center}
\includegraphics[height=6.0cm, width=5.5cm, scale=1]{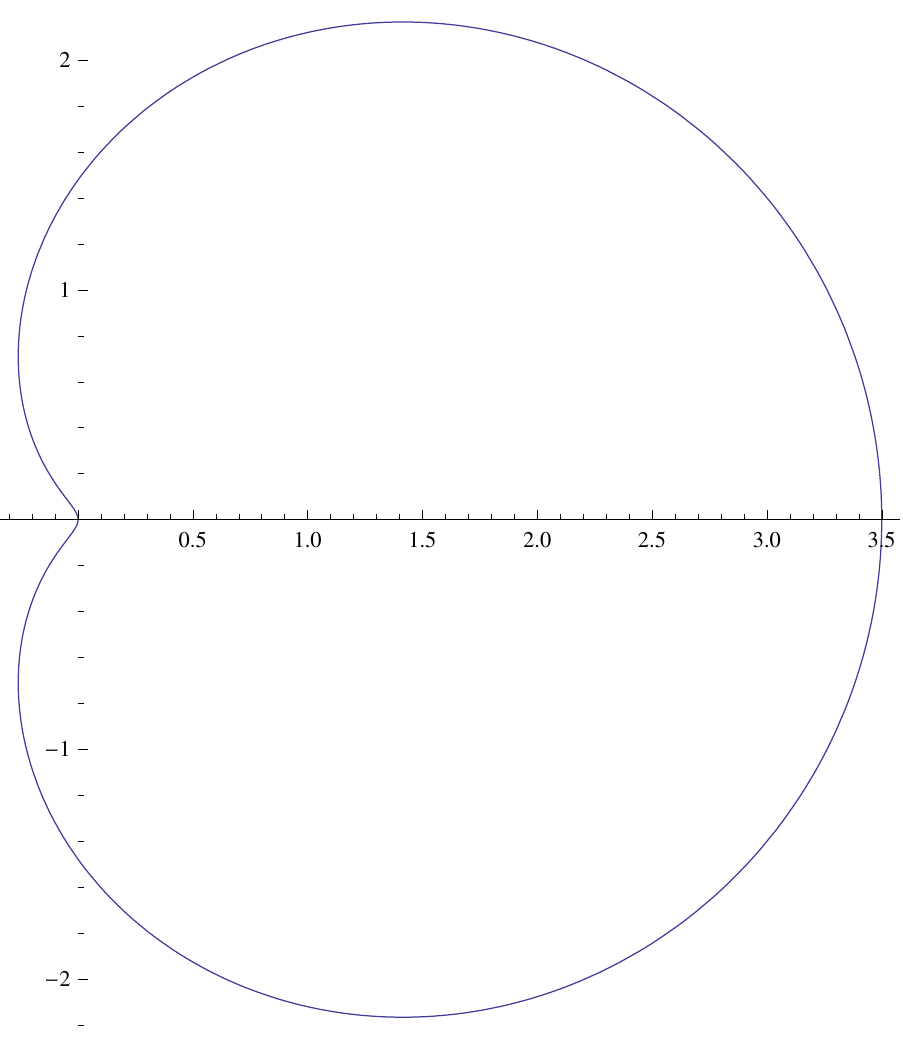}
\hspace{1cm}
\includegraphics[height=6.0cm, width=5.5cm, scale=1]{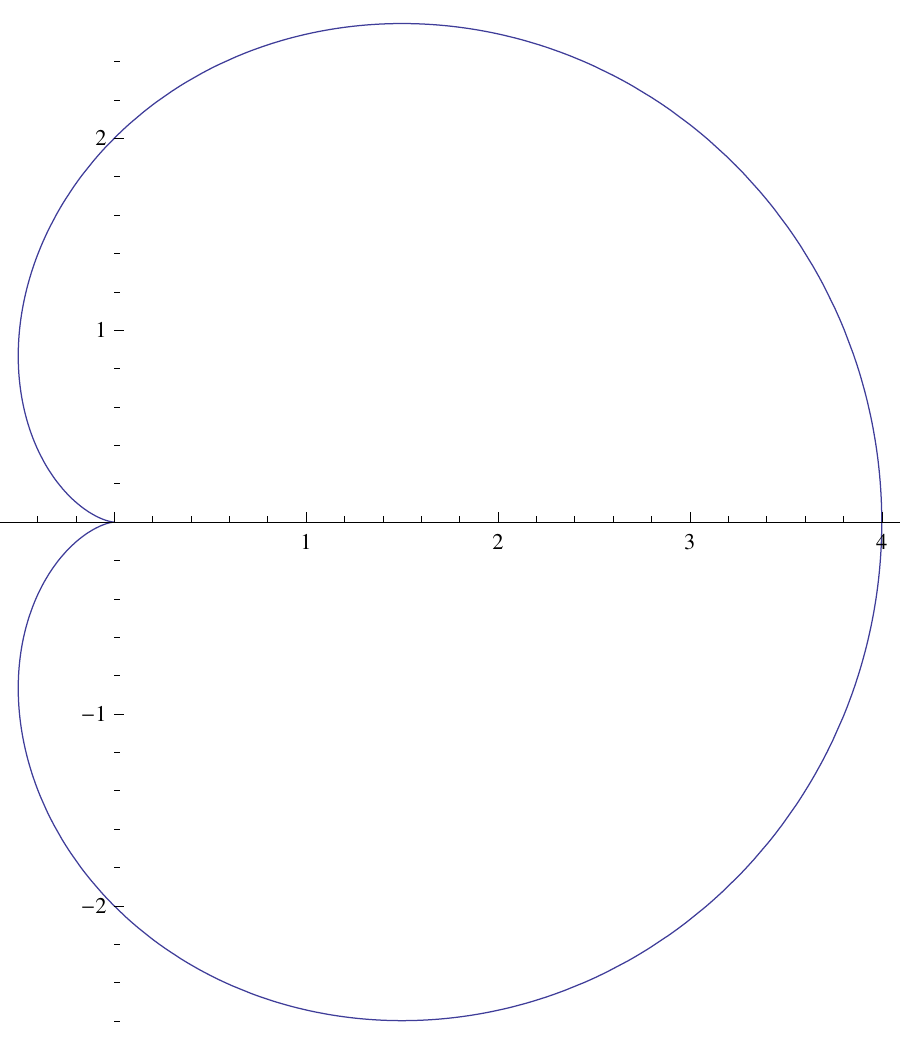}
\end{center}
$\lambda=0.75$ \hspace{5cm} $\lambda=1$
\caption{Graph of  $f(\lambda)=1+(1+\lambda)e^{i\alpha}+\lambda e^{2i\alpha}$ for certain values of $\lambda$, where $0\leq\alpha\leq2\pi$.\label{limacon2}}
\end{figure}

Now, we turn to our second duty. Since $s$ is injective in $\mathbb{D}$, we have to show that the image of $\mathbb{D}$ under the functions $z/g$ defined by (\ref{f1}) with $|\beta |\leq \beta_1$ is contained in the domain bounded by the  lima\c{c}on
\[
1+(1+\lambda)e^{i\alpha}+\lambda e^{2i\alpha},\quad \alpha \in [0,2\pi].\]
By calculations similar to the above ones, we see that this is equivalent to the assertion that for $|\beta |\leq \beta_1$ the points
\[\{lz+\lambda z^2:\, z\in \mathbb{D}\},
\]
are contained in the set
\[
\{w:\, w =e^{i\beta}((1+\lambda)u+\lambda u^2), ~u\in \mathbb{D}\}.
\]
This is a simple consequence of the fact that $(-1,0)$ is the point nearest to the origin of the  lima\c{c}on (see Figure \ref{limacon2})
\[
(1+\lambda)e^{i\alpha}+\lambda e^{2i\alpha},\quad \alpha \in [0,2\pi],\]
and that the point of intersection of this  lima\c{c}on turned around with angle $\beta_1$, the unit disk and the lima\c{c}on
\[
le^{i\alpha}+\lambda e^{2i\alpha},\quad \alpha \in [0,2\pi],\]
is the point $e^{-i\beta_1}.$ This completes the proof of \eqref{OPW7-eq4c}.

For the proof of the second part, by the definition of subordination, we simply rewrite \eqref{OPW7-eq4c} as
$$ \frac{z}{f(z)}= 1+ (1+\lambda) \omega (z)+\lambda \omega ^{2}(z),
$$
where $\omega$ is analytic in $\ID$ and $|\omega (z)|\leq |z|.$ It follows that
from the last equality that
$$\left|\frac{z}{f(z)}-1\right|\leq -1+1+(1+\lambda )|z|+\lambda |z|^{2} =-1+(1+\lambda |z|)(1+|z|)
$$
and the proof is complete.
\epf

According to Theorem \ref{OPW7-th4}, one has the estimate
$$\left|\frac{z}{f(z)}\right|\leq (1+\lambda r)(1+r)~\mbox{ for $|z|=r$}
$$
for $f\in {\mathcal U}(\lambda)$, $\lambda \in (0,1]$.

\br
We remark that Theorem \ref{OPW7-th2} follows from Theorem \ref{OPW7-th4}. Indeed,
there is nothing to prove if $\lambda =1$. Thus, if  $f\in {\mathcal U}(\lambda)$ for some $0<\lambda <1$, then we have
$$\frac{z}{f(z)}\prec 1+(1+\lambda)z+\lambda z^2.
$$
By Rogosinski's theorem \cite{Rogo43} (see also \cite[Theorem 6.2]{Duren:univ}), it follows that
$$1+ |a_2|^2\leq 1+ (1+\lambda)^2
$$
which implies that  $|a_2|\leq 1+\lambda$ for $\lambda \in (0,1)$.
\er

Under a mild restriction on $f$, one could improve the bound $|a_2|\leq 1+\lambda$ by establishing a region of
variability of $a_2$. In the next result we deal with this.

\bthm\label{OPW7-th2a}
Let $f\in {\mathcal U}(\lambda)$ for some $0<\lambda \leq 1$, and such that
\be\label{eq4c}
\frac{z}{f(z)}\neq (1-\lambda )(1+z),\,\,z\in \ID.
\ee
Then, we have
\be \label{eq4d}
\frac{z}{f(z)}-(1-\lambda )z\prec 1+2\lambda z+\lambda z^{2}
\ee
and the estimate $ |a_2-  (1-\lambda )|\leq 2\lambda $ holds. In particular, $ |a_2| \leq 1+\lambda $ and the estimate is sharp as the
function $f_{\lambda}(z) =z/((1+\lambda z)(1+z))$ shows.
\ethm
\bpf
Notice that there is nothing to prove if we allow $\lambda =1$.
Let $f \in \mathcal{U}(\lambda)$ for some $\lambda \in (0,1)$. Then, by the assumption \eqref{eq4c}, the function $g$ is defined by
\be\label{eq5c}
\frac{ z}{g(z)}=\frac{z}{f(z)}-\left(1-\lambda\right)(1+z),
\ee
has the property that $z/g(z)$ is non-vanishing and $g'(0)=1/\lambda$ and hence, it is easy to see that $G=\lambda g$ belongs to $\mathcal{U}$.
Consequently, by the last subordination relation in Theorem \ref{OPW7-th4}, we find that
$$ \frac{z}{G(z)}=
\frac{1}{\lambda}\left ( \frac{z}{f(z)}-(1-\lambda )(1+z) \right ) =1-\frac{a_{2}- (1-\lambda )}{\lambda } z+ \cdots \prec (1+z)^{2},
$$
which is obviously equivalent to \eqref{eq4d}.  The coefficient inequality $|(a_{2}- (1-\lambda ))/\lambda |\leq 2$ is a consequence of
Rogosinski's theorem. Thus, $|(a_{2}- (1-\lambda )|\leq 2\lambda$ holds.
\epf

It is not clear whether the condition \eqref{eq4c} is necessary for a function $f$ to belong to the family $\mathcal{U}(\lambda)$.

%\bcor\label{OPW7-cor3a}
%Suppose that $f(z)=z+\sum_{n=2}^{\infty}a_nz^n$ belongs to ${\mathcal U}(\lambda)$ for some $0<\lambda <1$. Then, we have
%$$|a_3-a_2^2|\leq \sqrt{1+2\lambda +2\lambda ^2 -|a_2|^2}.
%$$
%\ecor
%\bpf
%Using  Theorem \ref{OPW7-th4}, we see that
%$$\frac{z}{f(z)}= 1-a_2z-(a_3-a_2^2)z^2+\cdots \prec 1+(1+\lambda)z+\lambda z^2.
%$$
%Again, by Rogosinski's theorem \cite{Rogo43} (see also \cite[Theorem 6.2]{Duren:univ}), it follows that
%$$|a_2|^2+ |a_3-a_2^2|^2 \leq (1+\lambda)^2 +\lambda ^2
%$$
%and the desired result follows.
%\epf

\bthm\label{OPW7-cor3a}
Suppose that $f(z)=z+\sum_{n=2}^{\infty}a_nz^n$ belongs to ${\mathcal U}(\lambda)$ for some $0<\lambda \leq 1$. Then, we have the sharp estimate
$$|a_3-a_2^2|\leq \lambda .
$$
\ethm
\bpf
It is a simple exercise to see that
$$\left(\frac{z}{f(z)}\right)^2f'(z)= 1+(a_3-a_2^2)z^2+\cdots = 1+\lambda z^2\omega (z)
$$
where $\omega\in{\mathcal B}$, i.e. $\omega$ is analytic in $\ID$ such that $|\omega(z)|\leq 1$ for $z\in\ID$. Hence,
we must have $|a_3-a_2^2|\leq \lambda.$
Equality is attained if and only if $\omega(z)=e^{i\theta}$ for some $\theta \in [0,2\pi]$, i.~e. for functions
$f\in {\mathcal U}(\lambda)$ of  the form
\begin{equation}\label{f1a}
f(z)=\frac{z}{1-a_2z- \lambda e^{i\theta}z^2} .
\end{equation}
To get all extremal functions, we consider all those functions, where we may assume $a_2\geq 0$. The condition
$$1-a_2z- \lambda e^{i\theta}z^2 \neq 0
$$
is equivalent to this condition. It is clear that this is fulfilled if $a_2\leq 1-\lambda$. For $1-\lambda <a_2\leq 1+\lambda$ we get by a
reasoning similar to that used in the proof of Theorem \ref{OPW7-th4} that the condition is fulfilled if and only if
\begin{equation}\label{f2-c}
\cos \theta\,\leq \frac{(1-\lambda^2)^2\,-\,a_2^2(1+\lambda^2)}{2\lambda a_2^2}.\end{equation}
Hence, the extremal functions are those of the form \eqref{f1a}, where in addition \eqref{f2-c} is satisfied.
\epf

We observe that for $\lambda = 1$, the above inequality leads to the well-known estimate $|a_3-a_2^2|\leq 1$ which holds for $f\in{\mathcal S}$
and the Koebe function $k(z)=z/(1-z)^2$ gives the equality.

\section{Marx type implication for functions in $\mathcal{U}$\label{sec4}}

According to Theorem \ref{OPW7-th4},   one has
$${\rm Re}\,\sqrt{\frac{f(z)}{z}}>\frac{1}{2},\,\,z\in\ID,
$$
if $f\in \mathcal{U}$. This result is known to be true also for functions in the family $\mathcal{S}^{\star}$ of starlike functions in $\ID$
(see Marx \cite{Marx}) although the class  $\mathcal{U}$ neither contains  $\mathcal{S}^{\star}$ nor is contained in $\mathcal{S}^{\star}$.
On the other hand, since the structure of the class $\mathcal{U}$ allows us to determine the lower bound for the functional
${\rm Re}\,\sqrt{f(z)/z}$,   as a function of the second Taylor coefficient $a_2$, it is natural to solve the problem of
finding  $\alpha=\alpha (|a_2|)\geq 1/2$ such that $f\in \mathcal{U}$ implies that
$${\rm Re}\,\sqrt{\frac{f(z)}{z}}>\alpha,\,\,z\in\ID.
$$
In the next theorem, we present a solution to this problem. Also, in our result below, we observe that $\alpha (2)=1/2$
which is indeed the correct bound as the Koebe function $z/(1-z)^2$ shows. However, we could not claim that the bound $\alpha (|a_2|)$ is best possible.

\bthm\label{13-th7}
Let $f\in {\mathcal U}$  and $a_2 =f''(0)/2$. Then
$$ {\rm Re}\sqrt{\frac{f(z)}{z}}>\alpha( |a_{2}|)~ \mbox{ for $z\in\ID$},
$$
where
\be\label{eq8}
\alpha (x)=\frac{20+x-\sqrt{x^2+40x+16}}{24}, ~ 0\leq x\leq 2.
\ee
\ethm
\bpf
%The result is true for $a_2=0$ and hence, we assume that $a_2\neq 0$.
We recall from Lemma \ref{lem1} that the family  ${\mathcal U}$ is invariant under rotation and thus, it suffices to prove the theorem for
functions $f\in {\mathcal U}$ such that $a_{2}$ is real and non-negative and thus, throughout the proof, we may assume that $0\leq a_2\leq 2$.
Observe that $\alpha (x)$ is a decreasing function of $x\in [0,2]$ with $\alpha ([0,2]) = [1/2,2/3]$. We now let
\be\label{eq9}
\sqrt{\frac{f(z)}{z}} =p(z)= 1+\beta z+ \cdots ,
\ee
where $p$ is analytic in $\ID$,  $p(0)=1$ and $a_2$ is fixed and $0\leq \beta :=(a_2/2) \leq 1$.

We wish to  prove that
$$p(z)\prec q(z) := \frac{1+(1-2\alpha)z}{1-z}=1+2(1-\alpha)z+\cdots ,
$$
where $\alpha = \alpha(a_{2})$ is defined by \eqref{eq8}. We prove this by the method of contradiction.

Suppose that $p(z)$ is not subordinate to $q(z)$. Then, according to the result of  Miller and Mocanu \cite{MM-MMJ81,MM-MMJ85} (see also \cite{AliNagRavi2011}),
there exist points  $z_{0}\in \ID$ and $\zeta_{0}\in \partial (\ID)\backslash\{1\}$  such that
\be\label{eq10a}
p(z_{0})=q(\zeta_{0})~\mbox{ and }~ z_{0}p'(z_{0})=m\zeta_{0}q'(\zeta_{0}),
\ee
where
\be\label{eq10}
m\geq 1+\frac{q'(0)-\beta}{q'(0)+\beta} =\frac{8(1-\alpha)}{4(1-\alpha)+a_2}.
\ee
We notice that  $0\leq \beta=\frac{1}{2}a_{2}\leq q'(0)=2(1-\alpha)$. Also, we see that
\be\label{eq10b}
q(\zeta_{0})= \alpha +(1-\alpha)\frac{1+\zeta_0}{1-\zeta_0}=:\alpha +i\rho, \quad \rho \in \IR,
\ee
and a computation gives
\be\label{eq10c}
\zeta_{0}q'(\zeta_{0})= \frac{2(1-\alpha)\zeta_0}{(1-\zeta_0)^2}=-\frac{[(1-\alpha)^2+\rho^ 2]}{2(1-\alpha)}.
\ee
Further, using \eqref{eq9} and \eqref{OS-eq2}, it follows easily that
$$U_f(z)= \frac{1}{p^2(z)}+\frac{2zp'(z)}{p^3(z)}-1
$$
and thus, by \eqref{eq10a}, we obtain that
$$U_f(z_0)= \frac{1}{q^3(\zeta _0)} \left [q(\zeta _0) + 2m\zeta_{0}q'(\zeta_{0})-q^3(\zeta _0) \right ].
$$
By \eqref{eq10b} and \eqref{eq10c}, we deduce that
\beqq
|U_f(z_0)|^{2} &=& \frac{1}{|q(\zeta _0)|^6}\left|q(\zeta _0) + 2m\zeta_{0}q'(\zeta_{0})-q^3(\zeta _0)\right|^{2}\\
&=& \frac{1}{(\alpha^{2}+\rho^{2})^{3}}\left |\alpha +i\rho -\frac{m[(1-\alpha)^2+\rho^ 2]}{1-\alpha} -(\alpha +i\rho)^{3}\right |^2
\eeqq
and a calculation shows that $|U_f(z_0)|^{2}=\Phi(\rho^2),$ where
$$ \Phi(t)=\frac{(a+bt)^{2}+ct(d+t)^{2}}{(1-\alpha)^{2}(\alpha^{2}+t)^{3}},
$$
with
$$t=\rho^{2},\, a=(1-\alpha)^{2}(m-\alpha(1+\alpha)),\,b=m-3\alpha (1-\alpha),\, c=(1-\alpha)^{2},\, d=1-3\alpha^{2}.
$$
Clearly the proof will be completed if we can show that $\Phi(t)\geq 1$ for all $t\ge 0$ under the assumption on $\alpha=\alpha( a_{2})$
given by  \eqref{eq8}.  The inequality $\Phi(t)\geq 1$ is equivalent to
\be\label{eq11}
At^{2}+Bt +C\geq 0,
\ee
where $t\geq 0$,
$$A=b^{2}+2cd-3\alpha^{2}(1-\alpha)^{2}, \, B=2ab+cd^{2}-3\alpha^{4}(1-\alpha)^{2},\,C=a^{2}-\alpha^{6}(1-\alpha)^{2}.
$$
In order to prove the inequality \eqref{eq11}, it suffices to show that $A,B, C$ are non-negative for $\alpha \in [1/2, 2/3]$. We begin to observe by \eqref{eq10} that
\beqq
a-\alpha^{3}(1-\alpha)&=&(1-\alpha)^{2}(m-\alpha(1+\alpha))-\alpha^{3}(1-\alpha)\\
&\geq& (1-\alpha)^{2}\left(\frac{8(1-\alpha)}{4(1-\alpha)+a_2}-\alpha(1+\alpha)\right)
-\alpha^{3}(1-\alpha)=0,
\eeqq
provided
\be\label{eq12}
\frac{8(1-\alpha)}{4(1-\alpha)+a_2}-\alpha(1+\alpha)=\frac{\alpha^{3}}{1-\alpha},
\ee
which is the same as  $ 12\alpha ^2-\alpha (20+a_{2})+8=0$. Solving this equation gives the
solution $\alpha =\alpha( a_{2})$ expressed by \eqref{eq8}, and hence, $C\geq 0$.
It remains to show that $A\geq 0$, $B\geq 0 $ for $\alpha \in [1/2, 2/3]$.
The last inequality shows that $a\geq \alpha^{3}(1-\alpha)$ and
\beqq
b &=& m-3\alpha (1-\alpha)\\
&\geq &\frac{8(1-\alpha)}{4(1-\alpha)+a_2}-3\alpha(1-\alpha), ~\mbox{ by \eqref{eq10},}\\
&=&\frac{\alpha}{1-\alpha} -3\alpha(1-\alpha), ~\mbox{ by \eqref{eq12},}\\
&=&\frac{\alpha[1-3 (1-\alpha)^2]}{1-\alpha} >0  ~\mbox{ for $\alpha \in [1/2, 2/3]$.}
\eeqq
Using these facts, we can prove that $A\geq 0$ for $\alpha \in [1/2, 2/3]$. We now find that
\beqq
A &=&b^{2}+2cd-3\alpha^{2}(1-\alpha)^{2}\\
&\geq & \left (\frac{\alpha}{1-\alpha} -3\alpha(1-\alpha)\right )^2 +2(1-\alpha )^{2}(1-3\alpha ^2) -3\alpha^2(1-\alpha)^2 \\
&=&\frac{\alpha ^2}{(1-\alpha)^2}- 6\alpha ^2 +2(1-\alpha)^2\\
&=&\frac{(2\alpha-1)^2(2-\alpha^2)}{(1-\alpha) ^2}
\eeqq
which is non-negative for $\alpha \in [1/2, 2/3]$. Similarly, we have
\beqq
B &=&2ab+cd^{2}-3\alpha^{4}(1-\alpha)^{2}\\
&\geq & 2\alpha^3(1-\alpha)\left (\frac{\alpha}{1-\alpha} -3\alpha(1-\alpha)\right ) +(1-\alpha )^{2}(1-3\alpha ^2)^2 -3\alpha^4(1-\alpha)^2 \\
&=&(2\alpha-1)^2(1+2\alpha -\alpha^2)
\eeqq
which is again non-negative for $\alpha \in [1/2, 2/3]$.

Finally, we have shown that $\Phi(t)\geq 1$, i.e. $|U_f(z_0)|\ge 1$, which is a contradiction to $|U_f(z)|<1$ in $\ID$ and hence to the
assumption that $p$ is not  subordinate to $q$.  Hence, we must have $p(z)\prec q(z)$ in $\ID$ which is equivalent to the desired result.
\epf

\section{Applications of Elementary Transformations \label{sec5}}

 Because each $f\in {\mathcal U}$ is non-vanishing in $\ID\backslash\{0\}$, $z/f(z)$  can be written as
\be\label{OS-eq1}
\frac{z}{f(z)}=1+\sum_{k=1}^{\infty}b_kz^k, \quad z\in \ID.
\ee
One of the sufficient conditions for functions $f$ of this form to belong to the class
${\mathcal U}$ is that (see \cite{OP-01,obpo-2007a})
\be\label{eq7a}
\sum_{n=2}^{\infty}(n-1)|b_{n}|\leq 1.
\ee
%It follows \cite{FR-2006,OP-01,PV2005} that neither ${\mathcal U}$ is included in ${\mathcal S}^{\star}$ nor
%includes ${\mathcal S}^{\star}$. Here ${\mathcal S}^{\star}$ denotes the class of functions $f\in{\mathcal S}$ such that
%$f(\ID)$ is starlike with respect to the origin.

\bthm\label{th1}
Let $ f\in {\mathcal A}$ and
$$ \frac{z}{f(z)}=1+b_{1}z+\sum_{n=2}^{\infty}(-1)^{n}b_{n}z^{n},
$$
where $b_{n}\geq 0$ for $n\geq 2$. Then $f \in {\mathcal S}$ if and only if $\sum_{n=2}^{\infty}(n-1)b_{n}\leq 1.$
\ethm\bpf
For $ f\in {\mathcal S}$, by Lemma \ref{lem1}, we have that $g(z)=-f(-z)\in {\mathcal S}$. Since
$$\frac{z}{-f(-z)} =1-b_{1}z+\sum_{n=2}^{\infty}b_{n}z^{n},
$$
then by the characterization given in \cite{obpo-2009a} (see also the survey article \cite{ObPo-Jammu}), $g\in {\mathcal U}$ if and only if
$\sum_{n=2}^{\infty}(n-1)b_{n}\leq 1$ if and only if $g\in {\mathcal S}$. The desired conclusion follows.
\epf

\bprob
It will be interesting to find necessary and/or sufficient conditions (as in \cite{obpo-2009a})
for the function $f\in {\mathcal A}$ of the following form to be univalent in $\ID$:
$$ \frac{z}{f(z)}=1+b_{1}z+\sum_{n=2}^{\infty}(-1)^{n-1}b_{n}z^{n}
~\mbox{ or }~ \frac{z}{f(z)}=1+b_{1}z-\sum_{n=2}^{\infty}b_{n}z^{n},
$$
where $b_{n}\geq 0$ for $n\geq 2.$
\eprob

A function $f$ analytic in $\ID$ is called $n$-fold symmetric $(n=1, 2, \ldots)$ if
$$f(e^{i2\pi/n}z)=e^{i2\pi/n}f(z)\quad \mbox{ for $z\in \ID$}.
$$
In particular, every $f\in{\mathcal A} $ is $1$-fold symmetric and every odd $f$ is $2$-fold symmetric.
Every $n$-fold symmetric function $f(z)=z+\sum_{k=2}^{\infty}a_kz^k$ can be written as
$$f(z)=z+a_{n+1}z^{n+1}+a_{2n+1}z^{2n+1}+ \cdots .
$$
Properties of various geometric subclasses of $n$-fold symmetric functions from $\mathcal S$ have been investigated by many authors \cite{Go}.
We now investigate certain analogous problems associated with the class ${\mathcal U}$.

\bthm\label{th2}
Let $ f\in {\mathcal U}$ be given by \eqref{OS-eq1}. Then for each $n\geq 2$, the function $f_n(z)$ defined by
$$\frac{z}{f_n(z)}=1+\sum_{k=1}^{\infty}b_{nk}z^{nk}
$$
also belongs to the class $\mathcal U$, whenever $z/f_n(z) \neq 0$ in $\ID$. More generally, if $ f\in {\mathcal U}(\lambda )$ is given by
\eqref{OS-eq1}, then $ f_n\in {\mathcal U}(\lambda )$ whenever it is non-vanishing in $\ID$.
\ethm
\bpf
Let $ f\in {\mathcal U}$ with $\phi (z) =z/f(z)$. Then $\phi (z)$ is nonvanishing and analytic in $\ID$ and has the form
$$ \frac{z}{f(z)}=\phi (z) =1+\sum_{k=1}^{\infty}b_kz^k.
$$
%Moreover,
%\be\label{OS-eq2}
%U_f(z)=\left (\frac{z}{f(z)} \right )^{2}f'(z)-1=\frac{z}{f(z)} -z\left (\frac{z}{f(z)} \right )'-1, \quad z\in\ID.
%\ee
Now, we define $\Phi_n$ by  $\Phi_n(z)=z/f_n(z)$ and $\omega=e^{i2\pi/n}$. Then, $\{\omega ^k: \, k=1,2, \ldots, n\}$
is the set of all $n$ $n$-th roots of unity. It is a simple exercise to see that
$$\Phi_n(z):=\frac{1}{n} \sum_{k=1}^{n}\phi (\omega ^kz)=\frac{1}{n} \sum_{k=1}^{n}\frac{z}{\omega ^{-k}f(\omega ^kz)}=1+\sum_{k=1}^{\infty}b_{nk}z^{nk}.
$$
%and in power series form this is equivalent to
%$$ \Phi_n(z)=1+\sum_{k=1}^{\infty}b_{nk}z^{nk} .
%%~\mbox{ and }~
%% \phi (\omega ^kz)=\frac{z}{\omega ^{-k}f(\omega ^kz)}.
%$$
Since $ f\in {\mathcal U}$, by Lemma \ref{lem1}, for each $k$, the function $F_k(z)$ defined by $F_k(z)=\omega ^{-k}f(\omega ^kz)$
clearly belongs to the class ${\mathcal U}$. By calculation and the relation \eqref{OS-eq2}, it follows that
$$U_{f_n}(z)=\frac{1}{n} \sum_{k=1}^{n}U_{F_k}(z)=\frac{1}{n} \sum_{k=1}^{n} \left [\left (\frac{\omega ^kz}{f(\omega ^kz)} \right )^{2}f'(\omega ^kz)-1\right ]
$$
and thus, $|U_{f_n}(z)|<1$ in $\ID$ for each $n\geq 2$. The proof is complete.
\epf

From the proof of  the following corollary, we see that the non-vanishing condition $f_n(z) \neq 0$ in $\ID$ in the
above theorem can be dropped for the case $n=2$.

\bcor\label{cor1}
If $f\in {\mathcal U}$, then the odd function $f_2$ defined by
$$ \frac{z}{f_2(z)}=\frac{1}{2}\left(\frac{z}{f(z)}+\frac{z}{-f(-z)}\right)
%=1+\sum_{k=1}^{\infty}b_{2k}z^{2k},
$$
also belongs to the class $\mathcal U.$ More generally, if $ f\in {\mathcal U}(\lambda )$, then $ f_2\in {\mathcal U}(\lambda )$.
\ecor \bpf
Let $f\in {\mathcal U}$. Then, by Lemma \ref{lem1}, $F$ defined by $F(z)= -f(-z)$ belongs to
${\mathcal U}$. Moreover, the condition $f(z)-f(-z)\neq 0$ for $z\in \ID\setminus \{0\}$ is satisfied,
because if $f(z)=f(-z)$ for some $z\in \ID \setminus\{0\}$, then, since $f$
is univalent, we have $z=-z$, i.e. $z=0$, which is a contradiction. Consequently,
$$\frac{z}{f_2(z)}=\frac{z^2}{f(z)f(-z)}\left(\frac{f(z)-f(-z)}{2}\right)
$$
is non-vanishing in $\ID$. Moreover, a calculation gives that if
$f\in {\mathcal U}$ is given by  \eqref{OS-eq1}, then $f_2$ takes the form
$$\frac{z}{f_2(z)}=1+\sum_{k=1}^{\infty}b_{2k}z^{2k}
$$
and thus, by Theorem \ref{th2}, $f_2\in {\mathcal U}$.
\epf

From the proof of Theorem \ref{th2}, the following general result could be proved easily and so, we omit its details.

\bcor\label{th3}
Let $ g_k\in {\mathcal U}(\lambda _k)$ for $k=1,2,\ldots ,n$ and $\mu _k, \lambda _k\in [0,1]$ for $k=1,2,\ldots ,n$ such that
$\mu_1\lambda _1+\cdots +\mu_n\lambda _n=1$.
If $\Phi$ defined by
$$\Phi (z)= \sum_{k=1}^{n}\mu_k \frac{z}{g_k(z)} =\frac{z}{\Psi (z)}
$$
is non-vanishing in $\ID$, then the function $\Psi (z)=\frac{z}{\Phi (z)}$ belongs to the class $\mathcal U$.
\ecor
\bpf It suffices to observe that
$$U_{\Psi}(z)=\sum_{k=1}^{n} \mu _kU_{g_k}(z)
$$
and the rest follows by taking the modulus on both sides and use the triangle inequality.
\epf

\bcor
Let $f\in \mathcal {U}$ be given by  \eqref{OS-eq1}. For $\theta \in [0,2\pi)$, the functions $f_{3}$ and $f_{4}$ defined by
$$\frac{z}{f_{3}(z)}=1+\sum_{n=1}^{\infty}b_{n}\cos(n\theta)z^{n}
~\mbox{ and }~\frac{z}{f_{4}(z)}=1+\sum_{n=1}^{\infty}b_{n}\sin(n\theta)z^{n}
$$
also belong to the class $\mathcal {U}$ (whenever $z/f_{3}$ and $z/f_{4}$ are non-vanishing in $\ID$).
\ecor
\bpf Lemma \ref{lem1} shows that the functions
$g_1(z)=e^{-i\theta}f(ze^{i\theta})$ and $g_2(z)=e^{i\theta}f(ze^{-i\theta})$ belong to the class
$\mathcal {U}$ and so does its convex combination (by Corollary \ref{th3} with $\mu_1=\mu_2=1/2$ and $\lambda_1=\lambda_2=1$).
Moreover, it follows from the power series representation of
$z/f(z)$ that
$$\frac{z}{f_{3}(z)}=\frac{1}{2}\left(\frac{z}{e^{-i\theta}f(ze^{i\theta})}
+\frac{z}{e^{i\theta}f(ze^{-i\theta})}\right) =1+\sum_{n=1}^{\infty}b_{n}\cos (n\theta)z^{n}
$$
from which we conclude that $f_3\in\mathcal {U}$, by Corollary \ref{th3}.

In order to prove that $f_{4}$ belongs to $\mathcal {U}$, we first observe that
$$\frac{z}{f_{4}(z)}=1+ \frac{1}{2i}\left(\frac{ze^{i\theta}}{f(ze^{i\theta})}
-\frac{ze^{-i\theta}}{f(ze^{-i\theta})}\right) =1+\sum_{n=1}^{\infty}b_{n}\sin(n\theta)z^{n},
$$
and, by a computation, we have
$$\left|U_{f_{4}}(z)\right|
 =\left|\frac{1}{2i}\left( U_{f}(ze^{i\theta})-  U_{f}(ze^{-i\theta})\right)\right|
\leq \frac{1}{2}\left(|U_{f}(ze^{i\theta})|+|U_{f}(ze^{-i\theta})|\right)<1,
$$
showing  that $f_{4}\in\mathcal {U}.$
\epf

In particular, if we set $\theta=\pi/2$, then $f_3(z)$ and $f_4(z)$ take the forms
$$\frac{z}{f_{3}(z)}=1-b_{2}z^{2}+b_{4}z^{4}- \cdots ~ \mbox{ and }~\frac{z}{f_{4}(z)}=1+b_{1}z-b_{3}z^{3}+ \cdots ,
$$
respectively, and thus, the above corollary provides us with new functions from $\mathcal {U}$.

\bthm
Let $ f\in {\mathcal U}$ be given by  \eqref{OS-eq1}. Then the function $g$ defined by
$$\frac{z}{g(z)}=1+\sum_{k=1}^{\infty}{\rm Re\,}\{b_{k}\}z^{k},
$$
with $z/{g(z)}\neq 0$ in $\ID$, also belongs to the class $\mathcal U.$
More generally, if $ f\in {\mathcal U}(\lambda )$, then $g\in {\mathcal U}(\lambda )$.
\ethm
\bpf Let $ f\in {\mathcal U}$. Then, by Lemma \ref{lem1}, $h(z)=\overline{f(\overline{z}})$ belongs to
${\mathcal U}$. Now, we observe that
\beqq
\frac{z}{g(z)}
= \frac{1}{2}\left[\left(1+\sum_{k=1}^{\infty}b_{k}z^{k}\right)
 +\overline{\left(1+\sum_{k=1}^{\infty}b_{k}\overline{z}^{k}\right)}\right]
=  \frac{1}{2}\left(\frac{z}{f(z)}+\frac{z}{h(z)}\right)
\eeqq
%\beqq
%\frac{z}{g(z)}
%%&=&1+\sum_{n=1}^{\infty}{\rm Re\,}{b_{n}\}z^{n}=1+\sum_{n=1}^{\infty}\frac{b_{n}+\overline{b_{n}}}{2}z^{n}\\
%&=& \frac{1}{2}\left[\left(1+\sum_{k=1}^{\infty}b_{k}z^{k}\right)
% +\overline{\left(1+\sum_{k=1}^{\infty}b_{k}\overline{z}^{k}\right)}\right]\\
%&=&  \frac{1}{2}\left(\frac{z}{f(z)}+\frac{z}{\overline{f(\overline{z})}}\right)\\
%&=&  \frac{1}{2}\left(\frac{z}{f(z)}+\frac{z}{g(z)}\right)
%\eeqq
and thus, we easily have
$$U_{g}(z)=\frac{z}{g(z)} -z\left (\frac{z}{g(z)} \right )'-1=\frac{U_{f}(z) +U_{h}(z)}{2}.
$$
Clearly, the last relation implies that $g\in {\mathcal U}.$
\epf

\bthm
Let $f\in \mathcal {U}$ be given by  \eqref{OS-eq1}. Then the function $F$ defined by
\be\label{eq3}
\frac{z}{F(z)}= 1+\sum_{n=1}^{\infty}b_{2n}z^{n}
\ee
belongs to the class $\mathcal {U}.$ More generally, if $ f\in {\mathcal U}(\lambda )$ is given by
\eqref{OS-eq1}, then $ F\in {\mathcal U}(\lambda )$.
\ethm
\bpf If $f\in \mathcal {U}$, then   we have the representation %(see \eqref{eq6})
\be\label{eq4}
\frac{z}{f(z)}=1+b_{1}z+z\int_{0}^{z}\frac{\omega (t)}{t^{2}}\,dt, \quad b_1=-a_2,
\ee
where  $\omega\in{\mathcal B}_1$. Here ${\mathcal B}_1$ denotes the
class of functions $\omega$ analytic in $\ID$ such that $\omega(0)=\omega '(0)=0$ and $|\omega(z)|<1$ for $z\in\ID$.
If we put
$$ \omega_{1}(z)=\int_{0}^{z}\frac{\omega (t)}{t^{2}}\,dt,
$$
then $\omega _1$ is analytic in $\ID$, $\omega_{1}(0)=0$ and $|\omega_{1}(z)|\leq |z|$. Moreover,
$|\omega'_{1}(z)| =|\omega(z)/z^2|\leq 1$ for every $z\in \ID$. Consequently, for $f\in \mathcal {U}$ one has
\be\label{eq5}
\frac{z}{f(z)}=1+b_{1}z+z\omega_{1}(z).
\ee
and thus, the function $\Psi$ defined by
$$ \Psi (z)=\frac{1}{2}\left(\frac{z}{f(z)}+\frac{-z}{f(-z)}\right) =1+\frac{z}{2}\left (\omega_{1}(z)-\omega_{1}(-z)\right )
$$
is analytic in $\ID$ and $|\Psi (z)-1|<1$ for $z\in \ID$. Consequently, $\Psi (z)\neq 0$ in $\ID$,
$$\Psi (z)=1+\sum_{n=1}^{\infty}b_{2n}z^{2n}
$$
and observe that $F$ defined by
$$\frac{z}{F(z)}= \Psi (\sqrt{z})= 1 -zW(z):=1+ \frac{z}{2}\left ( \frac{\omega_{1}(\sqrt{z})}{\sqrt{z}}-\frac{\omega_{1}(-\sqrt{z})}{\sqrt{z}}\right )
$$
is analytic in $\ID$, where $W$ is analytic in $\ID$. Next, we observe that
$$U_F(z)=\frac{z}{F(z)} -z\left (\frac{z}{F(z)} \right )'-1=z^2W'(z)
$$
and, in view of the fact that $|\omega (z)|\leq |z|^2$ and $|\omega'_{1}(z)|=|\omega(z)/z^2| \leq 1$, we can easily see that
$|z^2W'(z)|<1$ in $\ID$, which means that $F\in \mathcal{U}$.
\epf

\section{Some radius problem \label{sec6}}
 When we say that  $f\in\mathcal{U}$ in $|z|<r$ it means that the inequality
$|U_{f}(z)| < 1$ holds in the subdisk $|z|<r$ of $\ID$, which is indeed same as saying that $r^{-1}f(rz)$
belongs to the class $\mathcal{U}$.

\bthm
Let $f\in \mathcal {S}$ and $f$ be given by  \eqref{OS-eq1}. Then the function $F$ defined by
$$\frac{z}{F(z)}=1+\sum_{n=1}^{\infty}b_{2n}z^{n}
$$
belongs to the class $\mathcal {U}$ at least in the disk $|z|<r_{0} =0.778387$
(implying  $F$ is univalent in $|z|<r_{0}$), where $r_{0}\in (0,1)$ is the root of the equation
\be\label{eq7}
\frac{r(1 - r^2)^2}{2}\log\left(\frac{1+r}{1-r}\right)-(4+r^4-7r^2)=0.
\ee
\ethm
\bpf Assume that $f\in \mathcal {S}$ and is given by  \eqref{OS-eq1}.
In order to show that $F\in\mathcal {U}$ in the disk $|z|<r_{0}$, we need to prove that the function $G$ defined by
$G(z)=r^{-1}F(rz)$ belongs to $\mathcal {U}$ in $\ID$ for each $0<r\leq r_0$. Thus, we begin to consider
the function $G$ defined by
$$ \frac{z}{G(z)}=1+\sum_{n=1}^{\infty}b_{2n}r^{n}z^{n},
$$
where $0<r\leq1$. To prove  $G\in\mathcal {U}$,  by \eqref{eq7a},  it suffices to show that
$$S=:\sum_{n=2}^{\infty}(n-1)|b_{2n}|r^{n} \leq 1
$$
for $0<r\leq r_0$. To do this, we need to recall first the following inequality, namely, for $f\in \mathcal{S}$, the
necessary coefficient inequality  (\cite[Theorem 11 on p.193 of Vol. 2]{Go})
$$\sum_{n=2}^\infty (n-1)|b_n|^2\leq 1.
$$
This in particular gives that
$\sum_{n=2}^\infty (2n-1)|b_{2n}|^2\leq 1.
$
Now, we find that
\beqq
S&=&\sum_{n=2}^{\infty}\sqrt{2n-1}|b_{2n}|\frac{(n-1)}{\sqrt{2n-1}}r^{n}\\
&\leq& \left(\sum_{n=2}^{\infty}(2n-1)|b_{2n}|^{2}\right)^{\frac{1}{2}}
\left(\sum_{n=2}^{\infty}\frac{(n-1)^{2}}{2n-1}r^{2n}\right)^{\frac{1}{2}}\\
&\leq& \left(\sum_{n=2}^{\infty}\frac{(n-1)^{2}}{2n-1}r^{2n}\right)^{\frac{1}{2}}.
\eeqq
By a computation we see that
%$$\sum_{n=2}^{\infty}\frac{(n-1)^{2}}{2n-1}r^{2n} <
%\frac{1}{3}r^{4}+\sum_{n=3}^{\infty}\frac{n-7/2}{2}r^{2n} =
%\frac{1}{3}r^{4}+ \frac{8r^6-3r^8}{10(1- r^2)^2}=
%\frac{r^4(10+4r^2+r^4)}{30(1- r^2)^2}
%$$
%and thus, $S\leq 1$ holds provided
%$$\frac{r^4(10+4r^2+r^4)}{30(1- r^2)^2}\leq 1,
%$$
%$$r^4(10+4r^2+r^4)-30(1- r^2)^2=0.$$
\beqq
\sum_{n=2}^{\infty}\frac{(n-1)^2}{2n-1}r^{2n}&=&\frac{1}{2}\sum_{n=2}^{\infty}\left(n - \frac{3}{2} + \frac{1}{2(2n-1)} \right)r^{2n}\\
&=&\frac{1}{2}\left(\frac{r^2}{(1-r^2)^2}-r^2\right)-\frac{3r^4}{4(1-r^2)}- \frac{r^2}{4}+\frac{r}{8}\log\left(\frac{1+r}{1-r}\right)\\
&=&\frac{r^2 (3 r^2 - 1)}{4 (1 - r^2)^2}+\frac{r}{8}\log\left(\frac{1+r}{1-r}\right)
\eeqq
and thus, $S\leq 1$ holds provided
$$\frac{r^2 (3 r^2 - 1)}{4 (1 - r^2)^2}+\frac{r}{8}\log\left(\frac{1+r}{1-r}\right) \leq 1,
$$
i.e. if $0<r\leq r_{0}=0.778387$, where $r_{0}$ is the root of the equation \eqref{eq7}.
%$$\frac{r(1 - r^2)^2}{2}\log\left(\frac{1+r}{1-r}\right)-(4+r^4-7r^2)=0.
%$$
It means that $F$ is in the class $\mathcal {U}$ in the disc $|z|<r_{0}.$
\epf

In \cite{obpo-2005}, as a corollary to a general result, it has been shown that $|z|<1/\sqrt{2}$ is
the largest disk centered at the origin such that every function in $\mathcal{S}$ is included in $\mathcal{U}$. More precisely (see also \cite{Su2012}),
$$\sup\left\{r>0: \,r^{-1} f(rz)\in\mathcal{U} \mbox{ for every } f\in\mathcal{S} \right\}=1/\sqrt{2}.
$$
We conclude the paper with the following conjecture.

\bcon
If $f\in {\mathcal U}(\lambda)$ for some $0<\lambda \leq 1$.  Then $|a_n|\leq \sum_{k=0}^{n-1}\lambda ^k$.
\econ

There is nothing to prove if $\lambda =1$. Also, we have verified the truth of the conjecture for $n=3$.

\subsection*{Acknowledgement}
%\noindent
%{\bf Acknowledgement.}
The work of the first author was supported by MNZZS Grant, No. ON174017, Serbia.
The second author is currently on leave from Indian Institute of Technology Madras, India.

\end{document}